\documentclass[11pt]{article}
\usepackage{amsmath,amsthm}
\usepackage{pstricks,pst-node,pst-tree}
\usepackage{latexsym}
\usepackage{amssymb,amscd}
\usepackage{graphicx}
\usepackage{url}
\usepackage[left=1.5in,top=1in,right=1in,nohead]{geometry}
\numberwithin{equation}{section}
\let \:=\colon
\let \beg=\begin
\let \mb=\mathbb
\let \mc= \mathcal
\let \ra=\rightarrow
\let \st=\stackrel
\let \sub=\subset
\let \Ga=\Gamma
\let \La=\Lambda
\let \Sig=\Sigma
\let \la=\lambda

\let \al=\alpha
\let \be=\beta
\let \bet=\beta
\let \fl=\flushleft
\let \lfl=\lfloor
\let \rfl=\rfloor
\let \fr=\frac
\let \ga=\gamma
\let \rh=\rho
\let \om=\omega

\let \ot=\otimes
\let \op=\oplus
\let \ov=\overline

\let \de=\delta
\let \De=\Delta
\let \pr=\prime
\let \wt=\widetilde
\let \hra=\hookrightarrow
\let \lra=\longrightarrow

\begin{document}

\title{\bf Effective divisors on $\ov{\mc{M}}_g$ associated to curves with exceptional secant planes}
\author{\small Ethan Cotterill}
\maketitle

\beg{abstract}
This paper is a sequel to \cite{C}, in which the author studies secant planes to linear series on a curve that is general in moduli. In that paper, the author proves that a general curve has no linear series with exceptional secant planes, in a very precise sense. Consequently, it makes sense to study effective divisors on $\ov{\mc{M}}_g$ associated to curves equipped with secant-exceptional linear series. Here we describe a strategy for computing the classes of those divisors. We pay special attention to the extremal case of $(2d-1)$-dimensional series with $d$-secant $(d-2)$-planes, which appears in the study of Hilbert schemes of points on surfaces. In that case, modulo a combinatorial conjecture, we obtain hypergeometric expressions for tautological coefficients that enable us to deduce the asymptotics in $d$ of our divisors' virtual slopes.
\end{abstract}

\section{Introduction: Brill--Noether theory and divisors on $\ov{\mc{M}}_g$}
Determining when an abstract curve $C$ comes equipped with a map to $\mb{P}^s$ of degree $m$ is central to curve theory. The Brill--Noether theorem asserts that when the Brill--Noether number $\rho(g,s,m)$ is nonnegative, $\rho$ computes the dimension of the space of series $g^s_m$ on a general curve $C$ of genus $g$, and that when $\rho$ is negative, there are no $g^s_m$'s on a general curve. The main qualitative result of \cite{C} is a Brill--Noether-type result for {\it pairs} of linear series. A recap is as follows.

We say that an $s$-dimensional linear series $g^s_m$ {\it has a $d$-secant $(d-r-1)$-plane} provided an inclusion
\beg{equation}\label{seriesinseries}
g^{s-d+r}_{m-d}+ p_1 + \dots + p_d \hra g^s_m
\end{equation}
exists. Geometrically, \eqref{seriesinseries} means that the image of the $g^s_m$ intersects a $(d-r-1)$-dimensional linear subspace of $\mb{P}^s$ in $d$-points; such a linear subspace is a ``$d$-secant $(d-r-1)$-plane".

Next, let
\[
\mu(d,r,s):=d-r(s+1-d+r).
\]
The invariant $\mu$ computes the expected dimension of the space of $d$-secant $(d-r-1)$-planes along a fixed $g^s_m$. We prove \cite[Thm. 1]{C}:

{\it
If $\rho=0$ and $\mu=-1$, then a general curve $C$ of genus $g$ admits no linear series $g^s_m$ with $d$-secant $(d-r-1)$-planes.}

In fact, Farkas \cite{Fa} has proved, more generally, that when $\rho+ \mu<0$, a general curve $C$ of genus $g$ admits no linear series $g^s_m$ with $d$-secant $(d-r-1)$-planes.

An immediate corollary of the theorem is that {\it whenever $\rho=0$ and $\mu=-1$, curves with $d$-secant $(d-r-1)$-planes sweep out a divisor inside $\ov{\mc{M}}_g$}.

A couple of words are in order regarding the linear series parameter $r$. To avoid trivialities, we must have
\[
1 \leq r \leq s.
\]
Each specialization of $r$ defines an infinite family of examples, indexed by the incidence parameter $d$. The two most natural choices are $r=1$ and $r=s$, and in this paper we focus on the former, which corresponds to the situation studied by Lehn \cite{Le} in the context of the Hilbert scheme of points on a surface. Note that the case $r=1$ is ``tautological" in that it corresponds to the situation in which the evaluation map
\beg{equation}\label{eval}
V \ra H^0(L/L(-p_1-\dots-p_d))
\end{equation}
corresponding to a given linear series $(L,V)$ fails to be {\it surjective} along a $d$-tuple of points $p_1, \dots, p_d$ on the curve in question. In general, $r$ is the corank of the map \eqref{eval}.

\section*{Acknowledgements}
This work is an enhanced version of the second part of my doctoral thesis, which was carried out under the supervision of Joe Harris. I thank Joe, and also Steve Kleiman at MIT, for countless valuable conversations. The paper was completed while the author was a postdoctoral fellow at Queen's University and at the MPIM-Bonn. The author thanks both institutions for their hospitality and the excellent working conditions.

\subsection{Roadmap}\label{roadmap}
The material following this introduction is organized as follows. In Section~\ref{families}, we attempt to solve the problem of computing the expected number of linear series with exceptional secant planes in a given one-parameter family by computing the number of exceptional series along judiciously-chosen ``test families". Our general secant plane formula reads
\beg{equation}\label{generalsecfmla}
N^{d-r-1}_d= P_{\al} \al+ P_{\be} \be+ P_{\ga} \ga+ P_c c + P_{\de_0} \de_0
\end{equation}
where $\al$, $\be$, $\ga$, $c$, and $\de_0$ are certain tautological invariants associated to the given family.
Five relations are needed to determine the tautological coefficients $P_{\al}, P_{\be}$, $P_{\ga}, P_c$, and $P_{\de_0}$. Whenever $r=1$ or $r=s$, we find four out of the five relations needed; in general, we conjecturally obtain four out of five relations, with a fourth relation hinging on a conjecture about secant planes to K3 surfaces (Conjecture 1, Section $2.2$), which we prove under a curvilinearity hypothesis in Theorem~\ref{valueofNK3}. 

In Section~\ref{r=1}, we specialize to the case $r=1$, where our results are strongest. The basic reason for this is the existence of a comparatively  simple generating function for the expected number $N_d$ of $d$-secant $(d-2)$-planes to a $g^{2d-2}_m$, namely \cite[Thm~4]{C}:
\beg{equation}\label{Ndformula} 
\sum_{d \geq 0} N_d(g,m) z^d= \biggl(\fr{2}{(1+4z)^{1/2}+1} \biggr)^{2g-2-m} \cdot (1+4z)^{\fr{g-1}{2}}.
\end{equation}
Our derivation of \eqref{Ndformula} in our earlier paper \cite{C} follows from the fact that for all $d \geq 2$, $N_d$ the degree of a corank-one degeneracy locus of a certain map of vector bundles over the $d$th Cartesian product of a curve. Applied to this map, Porteous' formula expresses the degeneracy locus as a determinant in the Chern classes of the vector bundles. A combinatorial analysis enables us to reinterpret this determinant graphically, in terms of weighted counts of subgraphs of the complete labeled graph on $d$ vertices. An extension of the same analysis yields a graphical interpretation of the missing tautological coefficient $P_{\ga}$. Nevertheless, deducing a generating function for the coefficients $P_{\ga}$ on the basis of combinatorics alone seems difficult.

In Section~\ref{genfunctions}, we use our generating function for $N_d$, the relations among tautological coefficients already obtained from evaluating the general secant plane formula \eqref{generalsecfmla} along special families, and a conjectural fifth relation obtained via multiple-point formulas in Section~\ref{mpt}, to (conjecturally) determine generating functions for the tautological coefficients $P$, whenever $r=1$. In Section~\ref{hypergeomfns}, we use the generating functions determined in Section~\ref{genfunctions} in order to realize each of the tautological coefficients $P$ as linear combinations of generalized hypergeometric functions. 

In Section~\ref{mpt}, we deduce a conjectural fifth relation among tautological secant plane divisor coefficients whenever $r=1$ or $r=s$, by calculating secant plane formulas in a variety of particular cases. Our computations are based on an application of Kleiman's multiple point formula \cite{Kl} to the projection of an incidence correspondence of curves and secant planes onto a Grassmann bundle of secant planes. In Section~\ref{gsmexamples}, we list secant plane formulas in a number of particular examples. Most of these enumerative formulas are new.

Section~\ref{lebarz} is devoted to a discussion of Le Barz's cycle-theoretic approach to secant planes, as detailed in \cite{Lb1,Lb2}. We recast his calculations in a slightly larger degree of generality, and obtain generalizations of a normal bundle formula \cite[Prop. 3bis]{Lb1} of his. Our main purpose, however, is to lay the foundation for a future adaptation of Le Barz's approach to the setting of one-parameter families.

In Section~\ref{dprecap}, we review Deepak Khosla's computation of the Gysin pushforward from $A^1(\mc{G}^s_m)$ to $A^1(\widetilde{M}_{g,1})$, where $\widetilde{M}_{g,1} \sub \ov{\mc{M}}_{g,1}$ is a partial compactification of the space of smooth marked curves of genus $g$. Indeed, the present paper should be viewed as a companion to \cite{Kh2}. Applying Khosla's result, we compute the coefficients $b_{\la}$ and $b_0$ associated to the Hodge class and the boundary class of irreducible nodal curves, respectively, of secant plane divisor classes on $\ov{\mc{M}}_g$. As a consequence, we deduce in Section~\ref{slopexamples} that the slope of secant plane divisors is computed by $\fr{b_{\la}}{b_0}$ whenever $r=1$ or $r=s$ and $g \leq 23$. We then specialize to the case $r=1$, and use our hypergeometric formulas for tautological coefficients to prove, in Theorem~\ref{nonemptysecant}, that secant plane divisors on $\ov{\mc{M}}_g$ are nonempty when $r=1$. The class of each secant plane divisor depends on the incidence degree, $d$, as well as a second parameter, $a$. 

In Section~\ref{sloper=1}, we determine explicit formulas for the slopes of secant plane divisors in the case $r=1$, for small values of $a$. We also determine the asymptotics of the slope function as $d$ approaches infinity, for arbitrary (fixed) values $a$. 

In Section~\ref{bdrycoeffs}, we compute the coefficients $b_1$ and $b_2$ (corresponding to boundary classes $\de_1$ and $\de_2$, respectively)) of secant plane divisors on $\ov{\mc{M}}_g$, as functions of $b_{\la}$ and $b_0$. Our Theorem~\ref{delta2coeff} states that the pullback of any secant plane divisor class $\mbox{Sec}$ under the map $j_2: \ov{\mc{M}}_{2,1} \ra \ov{\mc{M}}_g$ given by attaching marked genus-2 curves to a general ``broken flag" curve is supported along the locus of curves with marked Weierstrass points.

\section{One-parameter families of curves with linear series}\label{families}
Given a complete curve $B$, let $\pi: \mc{X} \ra B$ denote a flat family of curves over $B$ whose generic fiber is smooth, with some finite number of special fibers that are irreducible curves with nodes. We equip each fiber of $\pi$ with an $s$-dimensional series $g^s_m$. That is, $\mc{X}$ comes equipped with a line bundle $\mc{L}$, and on $B$ there is a vector bundle $\mc{V}$ of rank $(s+1)$, such that 
\[
\mc{V} \hra \pi_* \mc{L}.
\]
If $\mu=-1$, we expect finitely many fibers of $\pi$ to admit linear series with $d$-secant $(d-r-1)$-planes. We then ask for a formula for the number of such series, given in terms of tautological invariants associated with the family $\pi$.

Ziv Ran's work on Hilbert schemes of nodal curves \cite{R1,R2} shows that the number of $d$-secant $(d-r-1)$-planes is a function $N^{d-r-1}_d$ of tautological invariants of the family $\pi$, namely:
\[
\al:=\pi_*(c_1^2(\mc{L})), \be:=\pi_*(c_1(\mc{L}) \cdot \om), \ga:=\pi_*(\om^2), \de_0, \text{ and } c:=c_1(\mc{V})
\]
where $\om= c_1(\om_{\mc{X}/B})$ and where $\de_0$ denotes the locus of points $b \in B$ for which the corresponding fiber $\mc{X}_b$ is singular.

More to the point, $N^{d-r-1}_d$ computes the degree of a the rank-$(d-r)$ locus for a map of vector bundles
\beg{equation}\label{evalbundlemap}
\mc{V} \ra \mc{T}^d(\mc{L})
\end{equation}
over the Hilbert scheme $\mc{X}^{[d]}_B$ of degree-$d$ subschemes of fibers of $\pi$. Above a point with support $(p_1,\dots,p_d)$, the map \eqref{evalbundlemap} is precisely the evaluation map \eqref{eval} mentioned earlier. For any fixed choice of $s$, Porteous'  formula, in tandem with Ran's calculus for intersections of tautological classes on $\mc{X}^{[d]}_B$, enables us to deduce that
\vspace{-5pt}
\beg{equation}\label{basicformula}
N^{d-r-1}_d= P_{\al} \al+ P_{\be} \be+ P_{\ga} \ga+ P_c c + P_{\de_0} \de_0
\end{equation}
where the arguments $P$ are polynomials in $m$ and $g$ with coefficients in $\mb{Q}$.

\medskip
Because such a formula \eqref{basicformula} in tautological invariants exists, the problem of evaluating it reduces to producing sufficiently many relations among the coefficients $P$. In \cite{C}, we obtain three of these. Namely, we have:
\beg{enumerate}
\item $2m P_{\al}+ (2g-2) P_{\be}+ (s+1) P_c=0$. 

This relation arises because \eqref{basicformula} is invariant under the renormalization $\mc{L} \mapsto \mc{L} \otimes \pi^* \mc{O}(D)$ that trivializes $\mc{V}$.
\item $P_c= -A(d,g,m)$, where $A(d,g,m)$ is the number of $d$-secant $(d-r)$-planes to a general curve of degree $m$ and genus $g$ in $\mb{P}^{s+1}$ that intersect a general line. 

This relation arises from evaluating \eqref{basicformula} along the family of projections of a general curve of degree $m$ in $\mb{P}^{s+1}$ from points along a disjoint line. 
\item $(-2m-2g)P_{\al}+ (2-2g)P_{\be}+ (-m-1)P_c= (d+1)A^{\pr}(d,g,m)$, where $A^{\pr}(d,g,m)$ is the number of $(d+1)$-secant $(d-r)$-planes to a general curve of degree $(m+1)$ and genus $g$ in $\mb{P}^{s+1}$.

This relation arises from evaluating \eqref{basicformula} along the family of projections of a generic curve of degree $m+1$ in $\mb{P}^{s+1}$ from points along the curve.
\end{enumerate}

Note that the families corresponding to the second and third relations are isotrivial. 

\subsection{A non-isotrivial family arising from K3 surfaces}
Let $S$ denote a $K3$ surface in $\mb{P}^s$, such that
\[
\mbox{Pic }S= \mb{Z}H \oplus \mb{Z}[C].
\]
where $H$ is the class of a hyperplane section, while $C$ is a smooth, irreducible curve of genus $g$ such that $C \cdot H=m$. Such surfaces exist, for a dense set of $(d,m,s)$, by \cite[Thm. 1.1]{Kn2}. The base locus of a pencil of curves of class $[C]$ consists of $[C]^2=(2g-2)$ points. Accordingly, we set
\[
\mc{X}=\mbox{Bl}_{2g-2 \text{ pts}} S \text{ and } B= \mb{P}^1,
\]
and let $\pi: \mc{X} \ra B$ denote the projection associated to a particular choice of base point.
Clearly, $c_1(\mc{L})=H$. Likewise, the relative dualizing sheaf of our family is given by
\[
\om_{\mc{X}/\mb{P}^1} = \om_{\mc{X}} \otimes \pi^* \mc{O}_{\mb{P}^1}(2).
\]
Now let $f$ denote the class of a fiber of $\pi$, and let $E_i, 1 \leq i \leq 2g-2$, denote the classes of the exceptional divisors of the blow-up $\mc{X} \ra S$. Then
\[
\beg{split}
\om&= K_{\mc{X}}+ 2f \\
&= \sum_{i=1}^{2g-2} E_i + 2\biggl([C] - \sum_{i=1}^{2g-2} E_i\biggr) = 2[C] - \sum_{i=1}^{2g-2} E_i.
\end{split}
\]
Whence,
\[
\beg{split}
\ga&= 4[C]^2+ \sum_{i=1}^{2g-2} E_i^2= 6g-6, \al = H^2=2s-2, \text{ and } \\
\be&= 2[C] \cdot H= 2m.
\end{split}
\]
We compute $\de_0$ as follows. Let $\mb{C}^2 \subset H^0(\mc{O}_S(C))$ denote the two-dimensional subspace of sections defining our pencil. Let $\mc{X}_2$ denote the fiber product $\mc{X} \times_{\mb{P}^1} \mc{X}$, equipped with projections $\pi_1$ and $\pi_2$ onto each of its factors. Now let
\[
\mc{E}:= (\pi_1)_* (\pi_2^* \mc{O}_{\mc{X}}(C) \otimes \mc{O}_{\mc{X}_2}/\mc{O}_{\mc{X}_2}(-\De));
\]
over $p \in \mb{P}^1$, $\mc{E}_p$ comprises sections of $\mc{O}_S(C)$ modulo those vanishing to order $2$ at $p$. 

Note that the singular fibers of $\pi$ comprise the locus where the evaluation map
\[
\mb{C}^2 \otimes \mc{O}_S \st{\mbox{ev}}{\lra} \mc{E}
\]
fails to be surjective. It follows that $\de_0=c_2(\mc{E})$. On the other hand, it is not hard to see that there is an exact sequence
\[
0 \ra \mc{O}_S(C) \otimes \mc{T}_S^* \ra \mc{E} \ra \mc{O}_S(C) \ra 0;
\]
it follows that
\[
\beg{split}
c_t(\mc{E}) &= c_t(\mc{O}_S(C)) \cdot c_t(\mc{O}_S(C) \otimes \mc{T}_S^*) \\
&= (1+ t[C]) \cdot (1+ t(\al_1+ 2[C])+ t^2(\al_1 [C] + [C]^2+ \al_2))
\end{split}
\]
where $\al_i= c_i(\mc{T}_S^*)$. We deduce that
\[
\de_0= 2 \al_1 [C]+ 3[C]^2+ \al_2.
\]
Here $\al_1= c_1(K_S) = \sum_{i=1}^{2g-2} E_i$, while
\[
\al_2 = 24
\]
is the topological Euler characteristic of $S$. It follows that $\de_0= 6g+18$.

Finally, the vector bundle $\mc{V} \ra \mb{P}^1$ is trivial, since the $\mb{P}^s$ to which the fibers of $\mc{X} \ra \mb{P}^1$ map is fixed. So $c=0$.

Therefore, the third family yields the relation
\beg{equation}\label{K3relation}
(2s-2) P_{\al}+ 2m P_{\be}+ (6g-6) P_{\ga}+ (6g+18) P_{\de_0}= N_{K_3}
\end{equation}
where $N_{K_3}$ denotes the number of fibers of $\pi$ with exceptional secant plane behavior. 

\subsection{The value of $N_{K_3}$}

If $r=1$, then $\mu=-1$ forces $d=2s-1$ and $d-r-1=s-2$, so $S$ admits no $d$-secant $(d-2)$-planes, by \cite[Thm. 1.1]{Kn1}. It follows that $N_{K_3}=0$ when $r=1$. At the other extreme, if $r=s$ then the assumption that $\mu=-1$ forces
\[
d=2s-1 \text{ and } d-r-1=s-2.
\]
By B\'{e}zout's theorem, the degree-$(2s-2)$ surface $S$ admits no $(s-2)$-planes, so again we have $N_{K_3}=0$.

For a general choice of $(d,m,r,s)$, the value of $N_{K_3}$ is unclear. However, we conjecture that the following is true.
\beg{conj}
Let $S \sub \mb{P}^s$ be a $K3$ surface with Picard group
\[
\mbox{Pic }S= \mb{Z}L \oplus \mb{Z}\La
\]
where
\[
L^2= 2s-2, \La^2= 2g-2, \text{ and } \La \cdot L=m.
\]
If 
\beg{equation}\label{numerics}
\rh(g,s,m)=0 \text{ and } \mu(d,s,r)=-1,
\end{equation}
then $S$ admits no $d$-secant $(d-r-1)$-planes, except possibly when $m=2s$ and $g=s+1$.
\end{conj}

{\fl \bf NB:} The hypothesis that $\rho(g,s,m)=0$ implies that
\beg{equation}\label{rhozero}
m=s(a+1), \text{ and } g=(s+1)a,
\end{equation}
for some positive integer $a$. When $a=1$, i.e., when $m=2s$ and $g=s+1$, the curves of class $L$ on $S$ are canonical curves. As soon as any such curve admits a $d$-secant $(d-r-1)$-plane, it admits an $r$-dimensional family of such secant planes. Consequently, those canonical curves with $d$-secant $(d-r-1)$-planes comprise a locus of codimension at least 2. As a result, the case $a=1$ has no bearing on our determination (in Sections 4-6) of the classes of secant plane divisors on $\mc{G}^s_m$ or $\ov{\mc{M}}_g$.

{\fl The} remainder of the section will be devoted to a proof of the following result.
\beg{thm}\label{valueofNK3}
With the same hypotheses on $S$ as in Conjecture 1, $S$ admits no curvilinear $d$-secant $(d-r-1)$-planes.
\end{thm}

\beg{proof}
Assume, for the sake of argument, that $S$ admits a curvilinear $d$-secant $(d-r-1)$-plane. Let $Z \sub S$ denote the corresponding subscheme. Since $Z$ is curvilinear, $Z$ is contained in a smooth hyperplane section $Y$ of $S$, not necessarily unique. 

Note $H^1(S,L)=0$, because $L$ is globally generated. Whence, the exact sequence defining $Z$ in $S$
\[
0 \ra L \ot \mc{I}_{Z/X} \ra L \ra L \ot \mc{O}_Z \ra 0
\]
induces an exact sequence
\[
0 \ra H^0(S,L \ot \mc{I}_{Z/S}) \ra H^0(S,L) \st{\mbox{ev}}{\lra} H^0(S,\mc{O}_Z) \ra H^1(S, L \ot \mc{I}_{Z/S}) \ra 0
\]
in cohomology. Here $h^0(S,\mc{O}_Z)=d$, and $\mbox{rk}(\mbox{ev})=d-r$ because $Z$ determines a $d$-secant $(d-r-1)$-plane to $S$, by assumption. It follows that 
\beg{equation}\label{dimisr}
h^1(S, L \ot \mc{I}_{Z/S})=r.
\end{equation}
On the other hand, we clearly have 
\[
L \ot \mc{I}_{Y/S} \cong \mc{O}_S,
\]
while the adjunction theorem on $S$ implies $\mc{I}_{Z/Y}(K_S+L) \cong \mc{O}_Y(K_Y-Z)$, i.e., that
\[
L \ot \mc{I}_{Z/Y} \cong \mc{O}_Y(K_Y-Z).
\]
It follows that the exact sequence of (twisted) ideal sheaves
\[
0 \ra L \ot \mc{I}_{Y/S} \ra L \ot \mc{I}_{Z/S} \ra L \ot \mc{I}_{Z/Y} \ra 0
\]
induces an exact sequence
\[
\beg{split}
H^1(S,\mc{O}_S) &\ra H^1(S,L \ot \mc{I}_{Z/S}) \ra H^1(Y,\mc{O}_Y(K_Y-Z)) \ra H^2(S,\mc{O}_S) \\
&\ra H^2(S,L \ot \mc{I}_{Z/S})
\end{split}
\]
in cohomology. Here $H^1(S,\mc{O}_S)=0$, while 
\[
\beg{split}
&H^2(S,\mc{O}_S) \cong H^0(S,K_S)^{\vee} \cong H^0(S,\mc{O}_S)^{\vee} \cong \mb{C}, \\
&H^1(Y,\mc{O}_Y(K_Y-Z)) \cong H^0(Y,\mc{O}_Y(Z)), \text{ and} \\
&H^2(S,L \otimes \mc{I}_{Z/S}) \cong H^2(S,L) \cong H^0(S,-L)^{\vee} =0.
\end{split}
\]
By \eqref{dimisr}, it follows that $Z$ defines a $g^r_d$ with $\rho(g,r,d) =-1$ along the canonical curve $Y$.

Note that Lazarsfeld's theorem \cite[Lem. 1.3]{La} states that provided there are no multiple or reducible curves of class $L$ on $S$, the $g^r_d$ defined by $Z$ is Brill--Noether general, which yields the desired contradiction. More precisely, Lazarsfeld shows that provided a certain vector bundle $\mc{F}$ admits no nontrivial endomorphisms, there are no multiple or reducible curves of class $L$ on $S$. Further, as was pointed out in \cite{FKP}, to show that $\mc{F}$ admits no nontrivial endomorphisms, it suffices to show that on $X$ there is no decomposition
\beg{equation}\label{decomp}
L = M+N
\end{equation}
where $M$ and $N$ are effective and verify $h^0(M) \geq 2, h^0(N) \geq 2$.

To see why, recall that the argument of \cite[Lem. 1.3]{La} establishes that if $\mc{F}$ admits nontrivial endomorphisms, then $c_1(\mc{F}^*)=L$ decomposes nontrivially as a sum of effective classes $M+N$, where
\[
M=c_1(\widetilde{M}) \text{ and }N=c_1(\widetilde{N})
\]
for suitably chosen coherent sheaf quotients $\widetilde{M}$ and $\widetilde{N}$ of $\mc{F}^*$. But Lazarsfeld also shows that $\mc{F}^*$ is generated by its global sections, so $\mbox{det }\widetilde{M}$ and $\mbox{det }\widetilde{N}$ are also generated by their global sections (and are nontrivial); it follows that $h^0(M) \geq 2$ and $h^0(N) \geq 2$. 

To show that no decomposition \eqref{decomp} exists, we assume the opposite and argue for a contradiction. Note that if a decomposition \eqref{decomp} exists, then because $\mbox{det }\widetilde{M}$ and $\mbox{det }\widetilde{N}$ are generated by their global sections, $h^1(M)=h^1(N)=0$, and the Riemann-Roch formula yields
\[
h^0(M)= 2+ \fr{1}{2} M^2 \text{ and } h^0(N)= 2+  \fr{1}{2} N^2.
\]
Since $h^0(M) \geq 2, h^0(N) \geq 2$, we have
\beg{equation}\label{ineq1}
M^2 \geq 0 \text{ and } N^2 \geq 0.
\end{equation}
On the other hand, we also have
\beg{equation}\label{ineq2}
M \cdot \La \geq 0 \text{ and } N \cdot \La \geq 0.
\end{equation}
Now let
\[
M= \al L + \bet \La \text{ and } N= (1-\al) L -\bet \La.
\]
Then
\[
\beg{split}
M^2 &= (\al L + \bet \La)^2 \\
&= \al^2(2s-2)+ \bet^2(2g-2)+ 2 \al \bet m \geq 0,\\
N^2 &= ((1-\al) L -\bet \La)^2 \\
&= (1-\al)^2(2s-2)+ \bet^2(2g-2)- 2(1-\al)\bet m \geq 0,\\
M \cdot \La &= (\al L + \bet \La) \cdot \La \\
&= \al m+ \bet(2g-2) \geq 0, \text{ and}\\
N \cdot \La &= ((1-\al) L -\bet \La) \cdot \La= (1-\al)m -\bet(2g-2) \geq 0.
\end{split}
\]

Note that the last two inequalities combine to yield
\beg{equation}\label{combineq}
0 \leq \al m+ \bet(2g-2) \leq m.
\end{equation}
There are now two cases to consider, namely: $(\al>0, \bet<0)$, and $(\al<0, \bet>0)$. The argument is virtually identical in either case; we present it in the first case.

First, observe that \eqref{combineq} implies that
\beg{equation}\label{1i}
-\fr{\bet}{\al} (2g-2) \leq m \leq -\fr{\bet}{(\al-1)}(2g-2).
\end{equation}
Similarly, the inequality deduced from $M^2 \geq 0$ above implies that
\beg{equation}\label{2i}
m \leq -\fr{\al}{\bet}(s-1)-\fr{\bet}{\al}(g-1).
\end{equation}
Now let $x= -\fr{\bet}{\al}>0$. Then \eqref{2i} may be rewritten as
\[
(g-1) x^2 -mx+ (s-1) \geq 0.
\]
The left-hand side of \eqref{1i} forces
\[
\beg{split}
x \leq \fr{m -\sqrt{m^2-4(g-1)(s-1)}}{2g-2}, \text{ i.e.,}\\
-\bet \leq \biggl(\fr{m -\sqrt{m^2-4(g-1)(s-1)}}{2g-2}\biggr) \al.
\end{split}
\]
The right-hand side of \eqref{1i} now forces
\[
\beg{split}
m &\leq (m -\sqrt{m^2-4(g-1)(s-1)}) \fr{\al}{\al-1}, \text{ i.e.,} \\
1 - \fr{1}{\al} &\leq 1 - \fr{\sqrt{m^2-4(g-1)(s-1)}}{m}, \text{ i.e.,} \\
\al &\leq \fr{m}{\sqrt{m^2-4(g-1)(s-1)}}.
\end{split}
\]
Next, we apply \eqref{rhozero}, with $a \geq 2$. We deduce that $\al \leq 1$ necessarily, except when $a=2$, when $\al=2$ is also a possibility.

Similarly, if $(\al<0, \bet>0)$, we conclude that $-\al \leq 1$ except possibly when $a=2$, when $\al=-2$ is also a possibility.

We now analyze the possibilities that remain.
\beg{itemize}
\item {\bf If $\al=1$}, then the left-hand side of \eqref{combineq} yields $-\bet \leq \fr{m}{2g-2}= \fr{s(a+1)}{2(s+1)a-2}$, which forces $\bet=0$.
\item Similarly, {\bf if $\al=0$}, then the right-hand side of \eqref{combineq} yields $\bet=0$.
\item {\bf If $\al=-1$}, then the right-hand side of \eqref{combineq} yields $\bet \leq \fr{m}{g-1}= \fr{s(a+1)}{(s+1)a-1}$, so that $\bet=0 \text{ or }1$. Then \eqref{combineq} forces $\bet=1$. But then $M \cdot L= (-L+ \La) \cdot L= m- (2g-2) \geq 0$ forces $m \geq 2g-2$, which contradicts \eqref{rhozero}.
\item {\bf If $a=2$ and $\al=-2$}, then the right-hand side of \eqref{combineq} forces $\bet \leq 2$. So either $(\al,\bet)=(-2,1)$, or $(\al,\bet)=(-2,2)$. But the left-hand side of \eqref{combineq} precludes $(\al,\bet)=(-2,1)$, and the condition that $N^2 \geq 0$ precludes $(\al,\bet)=(-2,2)$.
\end{itemize}

We conclude immediately.
\end{proof}
In light of the preceding result, to prove Conjecture 1 it would suffice to show that the curvilinear locus is open in the space of $d$-secant $(d-r-1)$-planes to $S$. Perhaps the smoothability tests developed recently by Erman and Velasco \cite{EV} can help shed light on this problem of deformation theory.

\section{The case $r=1$}\label{r=1}
As mentioned in the Roadmap, our results are strongest when $r=1$, in which case our secant plane formula counts $d$-secant $(d-2)$-planes to $g^{2d-1}_m$'s. In this situation, the basic invariants $A(d,g,m)$ and $A^{\pr}(d,g,m)$ defined in Section 2 satisfy
\[
A(d,g,m)= A^{\pr}(d-1,g,m),
\]
because $d$-secant $(d-1)$-planes to a curve $C$ of degree $m$ and genus $g$ in $\mb{P}^{2d}$ that intersect a disjoint line $l$ are in bijection with $d$-secant $(d-2)$-planes to the image $\pi_l(C)$ of $C$ under the projection with center $l$. On the other hand, letting $N_d(g,m):=A(d,g,m)$, we have
\beg{equation}\label{genfunctionNd}
\sum_{d \geq 0} N_d(g,m) z^d= \biggl(\fr{2}{(1+4z)^{1/2}+1} \biggr)^{2g-2-m} \cdot (1+4z)^{\fr{g-1}{2}}.
\end{equation}
by \cite[Thm~4]{C}. Equivalently,
{\small
\beg{equation}\label{expformNd}
\sum_{d \geq 0} N_d(g,m) z^d= \exp \biggl({\sum_{n>0} \fr{(-1)^{n-1}}{n} \biggl[\binom{2n-1}{n-1}m+ \biggl(4^{n-1}- \binom{2n-1}{n-1}\biggr) (2g-2)\biggr] z^n}\biggr). 
\end{equation}
}

\subsection{Graph-theoretic interpretation of tautological coefficients}\label{graphtheoretic}
As explained in \cite[Sec. 3.3]{C}, the coefficients of $m$ and $(2g-2)$ of the sum inside the exponential \eqref{expformNd} represent weighted counts of connected subgraphs of the complete labeled graph $K_d$ on $d$ vertices; the exponential computes an analogous tally of (possibly) disconnected graphs.

{\fl We} begin by revisiting the calculation of the sum inside the exponential formula for $N_d$ in \eqref{expformNd}. Recall that for all $d \geq 2$, $N_d$ computes the number of $d$-tuples along a fixed, general embedded curve $X \sub \mb{P}V^*$ with polarizing line bundle $L$ for which the corresponding evaluation map
\beg{equation}\label{evalmap}
V \ra H^0(L/L(-p_1- \dots- p_d))
\end{equation}
has rank equal to $(d-1)$.

{\fl More precisely, the} Cartesian product $X^d$ comes equipped with a {\it secant bundle} $T^d(L)$ whose fiber over $(p_1,\dots, p_d)$ is precisely the right side of \eqref{evalmap}. Porteous' formula expresses the class in $H^*(X^d,\mb{Q})$ corresponding to the rank-$(d-1)$ locus of $V \ra T^d(L)$ as the determinant
\beg{equation}\label{det}
\left| \beg{array}{ccccc}
c_1 & c_2 & \cdots  & c_{d-1} & c_d \\
1 & c_1 & \cdots  & c_{d-2} & c_{d-1} \\
\cdots & \cdots & \cdots & \cdots & \cdots \\
0 & \cdots & 0 & 1 & c_1
\end{array} \right|
\end{equation}
where $c_i$ denotes the $i$th Chern class of the secant bundle $T^d(L)$ over $X^d$. Ran \cite{R1} showed that the Chern polynomial of $T^d(L)$ is given by
\beg{equation}\label{secantsplitting}
c_t(T^d(L))= (1+ l_1 t) \cdot (1+ (l_2-\De_2)t) \cdots (1+ (l_d-\De_d)t)
\end{equation}
where $l_i, 1 \leq i \leq d$ is the pullback of $c_1(L)$ along the $i$th projection $X^d \ra X$, and $\De_j, 2 \leq j \leq d$ is the diagonal class supported along
\[
\De_j= \{(x_1,\dots,x_d) \in X^d | x_i=x_j \text{ for some } i<j \}.
\]
Whence, abusively writing $c_i$ for the corresponding expression obtained by setting each $l_i$ equal to zero yields
\[
c_i= (-1)^i s_i(\De_1,\De_2,\dots, \De_d)
\]
where $s_i=s_i(y_1,\dots,y_d)$ is the $i$th elementary symmetric function in $d$ indeterminates $y_1, \dots, y_d$. Here $\De_1:=0$.

\medskip
There are two main ingredients involved in evaluating the determinantal formula \eqref{det}. The first is the basic symmetric function identity
\beg{equation}\label{matrixid}
\left| \beg{array}{ccccc}
s_1 & s_2 & \cdots & s_{n-1} & s_n \\
1 & s_1 & \cdots & s_{n-2} & s_{n-1} \\
\cdots & \cdots & \cdots & \cdots & \cdots \\
0 & \cdots & 0 & 1 & s_1
\end{array} \right| = \mathop{\sum_{i_1, \dots, i_{d} \geq 0}}_{i_1+ \cdots i_{d}=n} y_1^{i_1} \cdots y_d^{i_{d}}.
\end{equation}

The second is the fact that on $X^d$, we have
\beg{equation}\label{expanddiag}
l_j \cdot \De_{i,j}= p_i^* m \{ \mbox{pt}_X \},
\end{equation}
and
\[
\De_{i,j}^2= -p_i^* \om_X \cdot \De_{i,j}= -(2g-2)p_i^* \{\mbox{pt}_X\} \cdot \De_{i,j}
\]
for every choice of pairs of distinct positive integers $2 \leq i,j \leq d$. Here
\beg{equation}\label{diag}
\De_j= \sum_{i=1}^{j-1} \De_{i,j}
\end{equation}
for every $2 \leq j \leq d$, where $\De_{i,j}$ is the ``small" diagonal supported along $d$-tuples of points whose $i$th and $j$th coordinates agree.

{\fl The} latter fact implies that the class of the {\it $d$-secant $(d-2)$-plane locus} inside $X^d$ is a sum of degree-$d$ monomials in the diagonal summands $\De_{i,j}$ and the $l_k$. This sum, in turn, may be reinterpreted as a weighted sum of (possibly disconnected) subgraphs of the labeled complete graph $K_d$, whose vertices correspond to large diagonals $\De_j, 1 \leq j \leq d$.

Namely, writing each large diagonal $\De_j$ as a sum of small diagonals \eqref{expanddiag} and expanding the determinant \eqref{det} accordingly yields a sum of monomials in the line bundles $l_j$ and the small diagonals $\De_{i,j}$. Any such monomial $m$ may be encoded as a (multi)graph, whose edges $e_{ij}$, which may appear with multiplicity, are in bijection with the small diagonals $\De_{ij}$ dividing $m$. Additionally, the vertex $v_j$ is colored red if and only if $l_j$ divides $m$. (Since $X$ is a curve, we have $l_j^2=0$ for all $j$, so each vertex is colored at most once.) Typical  colored and uncolored examples in the case $d=8$ are pictured below.

\includegraphics[scale=.25]{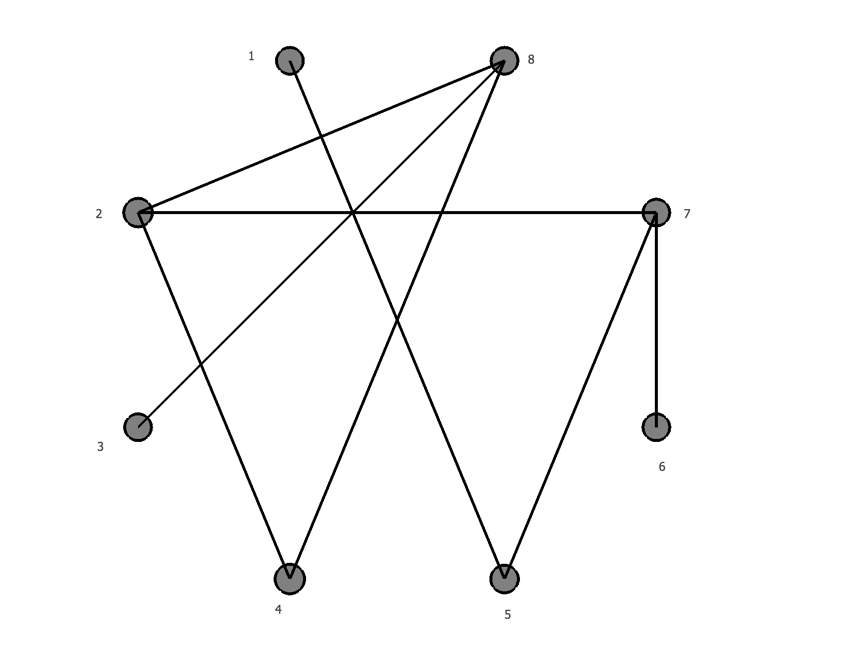}
\includegraphics[scale=.25]{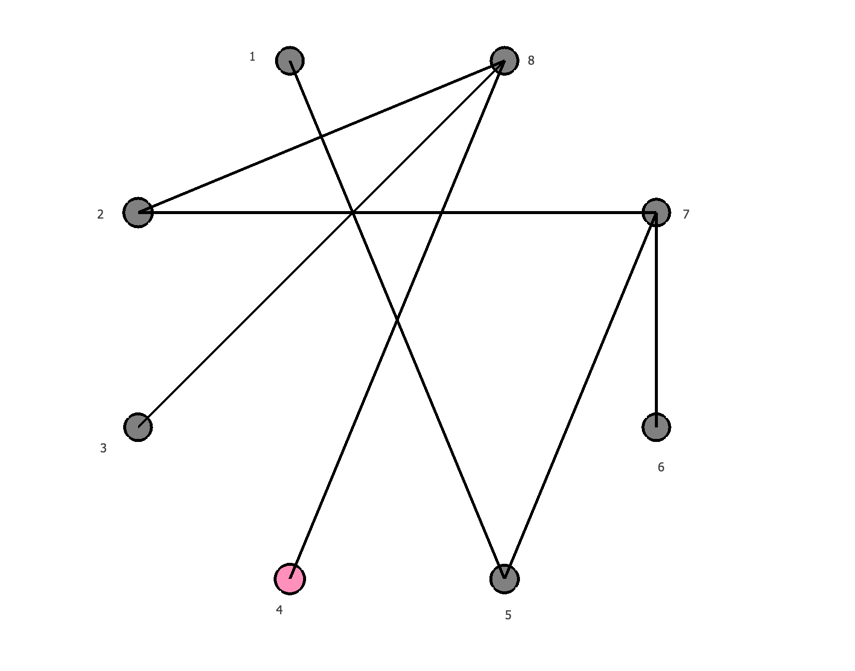}

The corresponding monomials are $\De_{1,5} \De_{2,7} \De_{2,4} \De_{2,8} \De_{3,8} \De_{4,8} \De_{5,7} \De_{6,7}$ and \\$l_4 \De_{1,5} \De_{2,7} \De_{2,8} \De_{3,8} \De_{4,8} \De_{5,7} \De_{6,7}$, respectively. 


{\fl Now} let $\mc{T}$ denote the set of connected spanning trees on $K_d$. To each vertex $v_j, 2 \leq j \leq d$ of a graph $G$ in $\mc{T}$, assign the weight
\[
w_{G,j}= (\mbox{indeg}(G,j))!.
\]
where $\mbox{indeg}(G,j)$ denotes the total indegree of the vertex $v_j$ in $G$. Now set $w_G= \prod_{2 \leq j \leq d} w_{G,j}$.
Then
\beg{equation}\label{Tweightedsum}
(2d-1) \sum_{G \in \mc{T}} w_G= \binom{2d-1}{d-1}(d-1)!
\end{equation}
corresponds to the coefficient of $m$ inside the exponential \eqref{expformNd}. 

{\fl Similarly}, let $\mc{S}$ denote the set of connected spanning graphs supported along $K_d$, with $d$ edges. To each vertex $v_j, 2 \leq j \leq d$ of a graph $G$ in $\mc{S}$, assign the following weight $w_{G,j}$:
\[
w_{G,j}= \binom{\mbox{indeg}(G,j)}{i_1,\dots, i_{j-1}}.
\]
Let $w_G= \prod_{2 \leq i \leq d} w_{G,j}$. Moreover, let $\mc{S}_1 \sub \mc{S}$ (resp. $\mc{S}_2 \sub \mc{S}$) comprise graphs all of whose edges appear with multiplicity $1$ (resp. graphs containing one edge appearing with multiplicity $2$). Clearly, $\mc{S}=\mc{S}_1 \cup \mc{S}_2$. Then
\beg{equation}\label{eqS1}
\sum_{G \in \mc{S}_1} w_G= \biggl(\sum_{i=0}^{d-3} \binom{2d-1}{i} \biggr) \cdot (d-1)!
\end{equation}
and
\beg{equation}\label{eqS2}
\sum_{G \in \mc{S}_2} w_G= \binom{2d-1}{d-2} \cdot (d-1)!.
\end{equation}
The total sum $\sum_{G \in \mc{S}_1} w_G+ \sum_{G \in \mc{S}_2} w_G$ corresponds to the coefficient of $(2g-2)$ inside the exponential \eqref{expformNd}.

{\fl The} Exponential Formula \cite[5.1.6]{St} implies that the generating function $\sum_{d>0} \wt{T}_d z^d$ for the aggregate count $\wt{T}_d$ of (possibly) disconnected subgraphs, each weighted by the number of times the corresponding monomial appears in the
expansion of the determinant \eqref{det}, is equal to the exponential of the corresponding generating function $\sum_{d>0} T_d z^d$ for {\it connected} graphs.


\subsection{Generating functions for tautological coefficients}\label{genfunctions}
Given our generating function \eqref{genfunctionNd} for $N_d(g,m)$ in tandem with the relations among tautological coefficients coming from our test families, determining generating functions for $P_{\al}=P_{\al}(d,g,m), P_{\be}=P_{\be}(d,g,m), P_{\ga}=P_{\ga}(d,g,m),$ and $P_{\de_0}=P_{\de_0}(d,g,m)$ is a purely formal matter. Namely, let
\beg{equation}\label{X,Y}
Z_{g,m}(z):=\biggl(\fr{2}{(1+4z)^{1/2}+1} \biggr)^{2g-2-m} \cdot (1+4z)^{\fr{g-1}{2}}.
\end{equation}
Then
\beg{equation}\label{explicitP}
\beg{split}
\sum_{d \geq 0} P_c(d,g,m) z^d&= Z_{g,m}(z), \\
\sum_{d \geq 0} P_{\al}(d,g,m) z^d&= Z_{g,m}(z) \biggl[ \fr{1}{2}- \fr{1}{2(1+4z)^{1/2}} \biggr], \text{ and}\\
\sum_{d \geq 0} P_{\be}(d,g,m) z^d &= Z_{g,m}(z) \biggl[ \fr{2z}{1+4z}- \fr{4z}{(1+4z)^{1/2} ((1+4z)^{1/2}+1)} \biggr].
\end{split}
\end{equation}
Moreover, we conjecturally also have
\beg{equation}\label{explicitPbis}
\beg{split}
\sum_{d \geq 0} P_{\ga}(d,g,m) z^d &=- \biggl(\fr{5}{6(2-g)}\biggr) \sum_{d \geq 0} P_{\al}(d,g,m) z^d + z \biggl(\fr{5}{6(2-g)}\biggr) \fr{d}{dz} \sum_{d \geq 0} P_{\al}(d,g,m)z^d \\
&+\biggl( \fr{(g+3)(m+3)}{12(2-g)}+ \fr{m}{12} \biggr) \sum_{d \geq 0} P_{\be}(d,g,m) z^d \\
&= Z_{g,m}(z) \biggl[\fr{z(32z^2-7(1+4z)^{3/2}+36z+7)}{6(1+4z)^{5/2}((1+4z)^{1/2}+1)} \biggr] \text{ and }\\
\sum_{d \geq 0} P_{\de_0}(d,g,m) z^d &=\biggl(\fr{1}{6(2-g)} \biggr) \sum_{d \geq 0} P_{\al}(d,g,m) z^d- z \biggl(\fr{1}{6(2-g)} \biggr) \sum_{d \geq 0} P_{\al}(d,g,m) z^d \\
&-\biggl(\fr{1}{12m}+\fr{(m-3)(g-1)}{12(2-g)} \biggr) \sum_{d \geq 0} P_{\be}(d,g,m) z^d \\
&= Z_{g,m}(z) \biggl[ \fr{z(32z^2-(1+4z)^{3/2}+12z+1)}{6(1+4z)^{5/2}((1+4z)^{1/2}+1)} \biggr].
\end{split}
\end{equation}

Our verification of the formulas \eqref{explicitPbis} for small values of $d$ is discussed in the next section. Now let
\[
X(z):= \fr{z(32z^2-7(1+4z)^{3/2}+36z+7)}{6(1+4z)^{5/2}((1+4z)^{1/2}+1)}, \text{ and } Y(z):= \fr{z(32z^2-(1+4z)^{3/2}+12z+1)}{6(1+4z)^{5/2}((1+4z)^{1/2}+1)}.
\]
The function $Y(z)$ has Taylor series
\[
\fr{1}{6}(3z^2 - 20 z^3  + 105 z^4  - 504 z^5  + 2310 z^6  - 10296 z^7  + 45045 z^8  - 194480 z^9+ \dots).
\]
In fact, it is not hard to show that
\[
[z^n] Y(z)= \fr{(-1)^{n-2}}{6} \cdot \fr{(2n-1)!}{n!(n-2)!}.
\]
Similarly, $X(z)$ has Taylor series
\[
\fr{1}{6}(-3z^2+ 28z^3- 177z^4+ 960z^5-4806z^6+ 22920z^7- 105837z^8+ 477688z^9- \dots),
\]
and we have
\[
[z^n] X(z)= (-1)^{n-1} \biggl(\fr{1}{2}\binom{2n}{n}(n+1)- 2^{2n-1}- \fr{1}{6} \cdot \fr{(2n-1)!}{n!(n-2)!} \biggr).
\]
Now set $\wt{g}=2g-2$. To conclude that $P_{\ga}$ and $P_{\de_0}$ have the generating functions predicted above, or equivalently, that the tautological coefficients satisfy the conjectural relation \eqref{finalrelation2}, it suffices to show the following.
\beg{conj}\label{pgaconst}
The exponential generating function for the constant terms of $\sum_{d \geq 0} P_{\ga}(d,\wt{g},m) z^d$ (resp., $\sum_{d \geq 0} P_{\de_0}(d,\wt{g},m) z^d$) is $X(z)$ (resp., $Y(z)$).
\end{conj}

To see why, note that Porteous' formula implies that $N_d^{d-2}$ is equal to the degree of the $(d+1) \times (d+1)$ determinant
\beg{equation}\label{det2}
\left| \beg{array}{ccccc}
c_1 & c_2 & \cdots  & c_{d} & c_{d+1} \\
1 & c_1 & \cdots  & c_{d-1} & c_{d} \\
\cdots & \cdots & \cdots & \cdots & \cdots \\
0 & \cdots & 0 & 1 & c_1
\end{array} \right|
\end{equation}
where $c_i$ denotes the $i$th Chern class of the secant bundle $\mc{T}^d(L)$ over $\mc{X}^{[d]}_B$. In fact, since the aim of this subsection is to calculate $P_{\ga}$, it suffices to consider families $\pi: \mc{X} \ra B$ all of whose fibers are {\it smooth}. Accordingly, we work on the relative Cartesian product $\mc{X}^d_B$, relative than the relative Hilbert scheme $\mc{X}^d_B$. The Chern polynomial of the (pullback) secant bundle now splits linearly exactly as in \eqref{secantsplitting}, except that now $l_i, 1 \leq i \leq d$ denote the pullbacks of $c_1(\mc{L})$ along the $i$th projection $\mc{X}^d_B \ra \mc{X}$, and $\De_j, 2 \leq j \leq d$ denote {\it relative} diagonal classes. In particular, the Chern classes of $\mc{T}^d(L)$ are, modulo $l_i$'s, elementary symmetric functions in the relative diagonals. Note that each relative diagonal $\De_j$ decomposes as a sum of irreducible diagonal classes
\[
\De_j= \sum_{i=1}^{j-1} \De_{i,j}
\]
supported along the loci of tuples $(x_1,\dots,x_d)$ in fibers of $\pi: \mc{X} \ra B$ for which $x_i=x_j$. We have the standard self-intersection formula
\beg{equation}\label{self-intersection}
\De_{i,j}^2= -\De_{i,j} \cdot \om_i
\end{equation}
where $\om_i$ denotes the $i$th pullback of the class $\om$ on $\mc{X}$. The symmetric function identity \eqref{matrixid} implies, that modulo $l_i$'s, the determinant \eqref{det2} computes $(-1)^{d+1}$ times
\[
\sum_{(i_1, \dots, i_{d-1} \geq 0} \De_2^{i_1} \cdot \De_d^{i_{d-1}}
\]
where the sum is over all degree-$(d+1)$ monomials in the relative diagonals.

Because $P_{\ga}$ computes the coefficient of $\om^2$, it follows that for all $d \geq 2$, $P_{\ga}(d,\wt{g},m)$ computes a weighted count of $(d+1)$-edged subgraphs $\mc{G}$ of the complete labeled graph on $d$ vertices with the following property.  Namely, we require that $\mc{G}$ possess a {\it unique} connected $(\widehat{d}+1)$-edged subgraph on $\widehat{d}$ vertices, for some $\widehat{d} \leq d$, and that the remaining connected components of $\mc{G}$ have either $d$ edges, or $(d-1)$ edges and a marked vertex corresponding to some pullback of $c_1(\mc{L})$. On the other hand, the coefficients of the exponential generating function for the constant terms of $P_{\ga}(d,\wt{g},m)$ compute weighted counts of $(d+1)$-edged {\it connected} subgraphs of the complete labeled graph on $d$ vertices. It follows from \cite[Prop. 5.1.1]{St} that $\fr{\sum_{d \geq 0} P_{\ga}(d,\wt{g},m) z^d}{Z_{\wt{g},m}(z)}$ is the exponential generating function for the constant terms of $P_{\ga}(d,\wt{g},m)$; the corresponding assertion about $\fr{\sum_{d \geq 0} P_{\de_0}(d,\wt{g},m) z^d}{Z_{\wt{g},m}(z)}$ is proved similarly.

\subsection{Hypergeometric formulas for tautological coefficients}\label{hypergeomfns}
Using the results of the preceding subsection, it is possible to realize $P_c, P_{\al}, P_{\be}$, $P_{\ga}$, and $P_{\de_0}$ as linear combinations of generalized hypergeometric series. Namely, we have the following result.

\beg{thm}\label{hypergeom} When $r=1$, the tautological secant-plane divisor coefficients $P_{\al}=P_{\al}(d,g,m)$, $P_{\be}=P_{\be}(d,g,m)$, and $P_c=P_c(d,g,m)$ are given by
\[
P_c=-\fr{g! (2g-2-m)!}{(g-2d)!d!(2g-2-m+d)!} { }_3F_2 \biggl[\beg{array}{ccc} -\fr{g}{2}+ \fr{m}{2}+ 1-d, & -\fr{g}{2}+ \fr{m+3}{2}-d, & -d \\ \fr{g+1}{2}-d, & \fr{g}{2}+1-d \end{array} \biggl| 1\biggr],
\]
\[
P_{\al}=\fr{g! (2g-2-m)!}{2(g-2d)!d!(2g-2-m+d)!} { }_3F_2 \biggl[\beg{array}{ccc} -\fr{g}{2}+ \fr{m}{2}+ 1-d, & -\fr{g}{2}+ \fr{m+3}{2}-d, & -d \\ \fr{g+1}{2}-d, & \fr{g}{2}+1-d \end{array} \biggl| 1\biggr]
\]
\vspace{-3pt}
\[
-\fr{(g-1)!(2g-2-m)!}{2(g-2d-1)!d!(2g-2-m+d)!} { }_3F_2 \biggl[\beg{array}{ccc} -\fr{g}{2}+ \fr{m}{2}+ 1-d, & -\fr{g}{2}+ \fr{m+1}{2}-d, & -d \\ \fr{g+1}{2}-d, & \fr{g}{2}-d \end{array} \biggl| 1\biggr],
\]
\[
P_{\be}= \fr{2(g-2)!(2g-2-m)!}{(g-2d)!(d-1)!(2g-3-m+d)!} { }_3F_2 \biggl[\beg{array}{ccc} -\fr{g}{2}+ \fr{m}{2}+ 1-d, & -\fr{g}{2}+ \fr{m+3}{2}-d, & 1-d \\ \fr{g+1}{2}-d, & \fr{g}{2}+1-d \end{array} \biggl| 1\biggr]
\]
\vspace{-3pt}
\[
-\fr{2(g-1)!(2g-1-m)!}{(g+1-2d)!(d-1)!(2g-2-m+d)!} { }_3F_2 \biggl[\beg{array}{ccc} -\fr{g}{2}+ \fr{m}{2}+ 1-d, & -\fr{g}{2}+ \fr{m+3}{2}-d, & 1-d \\ \fr{g}{2}+1-d, & \fr{g+3}{2}-d \end{array} \biggl| 1\biggr],
\]
Moreover, when $d \geq 3$, we have, assuming Conjecture 3:
\[
P_{\ga}= \fr{8(g-5)!(2g-1-m)!}{3(g+1-2d)!(d-3)!(2g-m+d-4)!} { }_3F_2 \biggl[\beg{array}{ccc} -\fr{g}{2}+ \fr{m}{2}+ 1-d, & -\fr{g}{2}+ \fr{m+3}{2}-d, & 3-d \\ \fr{g}{2}+ 1-d, & \fr{g+3}{2}-d \end{array} \biggl| 1 \biggr]
\]
\vspace{-3pt}
\[
-\fr{7(g-2)!(2g-1-m)!}{12(g-2d)!(d-1)!(2g-2-m+d)!} { }_3F_2 \biggl[\beg{array}{ccc} -\fr{g}{2}+ \fr{m}{2}+ 1-d, & -\fr{g}{2}+ \fr{m+1}{2}-d, & 1-d \\ \fr{g+1}{2}-d, & \fr{g}{2}+1-d \end{array} \biggl| 1 \biggr]
\]
\vspace{-3pt}
\[
+ \fr{3(g-5)!(2g-1-m)!}{(g-1-2d)!(d-2)!(2g-3-m+d)!} { }_3F_2 \biggl[\beg{array}{ccc} -\fr{g}{2}+ \fr{m}{2}-d, & -\fr{g}{2}+ \fr{m+1}{2}-d, & 2-d \\ \fr{g+1}{2}-d, & \fr{g}{2}-d \end{array} \biggl| 1 \biggr]
\]
\vspace{-3pt}
\[
+ \fr{7(g-5)!(2g-1-m)!}{12(g-3-2d)!(d-1)!(2g-2-m+d)!} { }_3F_2 \biggl[\beg{array}{ccc} -\fr{g+1}{2}+ \fr{m}{2}-d, & -\fr{g}{2}-1+ \fr{m}{2}-d, & 1-d \\ \fr{g}{2}-1-d, & \fr{g-1}{2}-d \end{array} \biggl| 1 \biggr],
\]
and
\[
P_{\de_0}= \fr{8(g-5)!(2g-1-m)!}{3(g+1-2d)!(d-3)!(2g-m+d-4)!} { }_3F_2 \biggl[\beg{array}{ccc} -\fr{g}{2}+ \fr{m}{2}+ 1-d, & -\fr{g}{2}+ \fr{m+3}{2}-d, & 3-d \\ \fr{g}{2}+ 1-d, & \fr{g+3}{2}-d \end{array} \biggl| 1 \biggr]
\]
\vspace{-3pt}
\[
-\fr{(g-2)!(2g-1-m)!}{12(g-2d)!(d-1)!(2g-2-m+d)!} { }_3F_2 \biggl[\beg{array}{ccc} -\fr{g}{2}+ \fr{m}{2}+ 1-d, & -\fr{g}{2}+ \fr{m+1}{2}-d, & 1-d \\ \fr{g+1}{2}-d, & \fr{g}{2}+1-d \end{array} \biggl| 1 \biggr]
\]
\vspace{-3pt}
\[
+ \fr{(g-5)!(2g-1-m)!}{(g-1-2d)!(d-2)!(2g-3-m+d)!} { }_3F_2 \biggl[\beg{array}{ccc} -\fr{g}{2}+ \fr{m}{2}-d, & -\fr{g}{2}+ \fr{m+1}{2}-d, & 2-d \\ \fr{g+1}{2}-d, & \fr{g}{2}-d \end{array} \biggl| 1 \biggr]
\]
\vspace{-3pt}
\[
+ \fr{(g-5)!(2g-1-m)!}{12(g-3-2d)!(d-1)!(2g-2-m+d)!} { }_3F_2 \biggl[\beg{array}{ccc} -\fr{g+1}{2}+ \fr{m}{2}-d, & -\fr{g}{2}-1+ \fr{m}{2}-d, & 1-d \\ \fr{g}{2}-1-d, & \fr{g-1}{2}-d \end{array} \biggl| 1 \biggr].
\]
The same formulas for $P_{\ga}$ and $P_{\de_0}$ hold when $d=2$, except that in each case the first hypergeometric summand should be suppressed. Finally, $P_{\ga}(1,g,m)=P_{\de_0}(1,g,m)=0$.
\end{thm}

\beg{proof}
See the proofs of \cite[Thm 4]{C} or \cite[Thm 5]{C2}.
\end{proof}

\subsection{Calculating the constant terms of $P_{\ga}(d,\wt{g},m)$}. We have verified Conjecture~\ref{pgaconst} for all $d \leq 8$, using an application of Kleiman's multiple-point formula. This approach is detailed in the next section. In this subsection, we describe a combinatorial approach which extends the graphical analysis of the preceding section. While interesting from a theoretical point of view, it involves a rather intricate recursion which makes it cumbersome to apply in practice. 


{\fl To begin}, write
\[
X(z)= \sum_{d \geq 0} \fr{(-1)^{d+1}}{d!} X_d z^d
\]
where $X_d$ is the weighted aggregate count of connected $(d+1)$-edged subgraphs of $K_d$, where each edge may appear with multiplicity at most 3. These subgraphs are classified naturally according to the multiplicities of their edges. Accordingly, $\mc{S}^{\pr}_i, i=2,3$ will designate those subgraphs that contain exactly {\it one} edge of multiplicity $i>1$. $\mc{S}^{\pr}_1$ will designate subgraphs without multiple edges, and $\mc{S}^{\pr}_{2,2}$ will designate subgraphs containing exactly {\it two} edges with multiplicity 2. 

{\fl We have}
\[
X_d= \sum_{G \in \mc{S}^{\pr}} w_G
\]
where $\mc{S}^{\pr}= \mc{S}^{\pr}_1 \cup \mc{S}^{\pr}_2 \cup \mc{S}^{\pr}_3 \cup \mc{S}^{\pr}_{2,2}$, and the weight $w_G$ is defined precisely as before.

We now analyze the contributions of subgraphs of each of these four types. Associated with each choice of partition $\la= (\la_1, \dots, \la_l)$ of $d+1$, there is a certain (unweighted) number $C^{\mc{S}^{\pr}}_{\la}$ of graphs $G \in \mc{S}^{\pr}$ with indegree sequence  $\lambda$ relative to some subset of $l$ vertices of $K_d$. For any choice of $G$, let $G^*$ denote the subgraph of $K_d$ obtained by deleting every edge of $G$ incident to $v_d$, let $i_G$ denote the number of these incident edges, and let $\ga_G$ denote the number of connected components of $G^*$. A useful observation for what follows is that $i_G-\ga_G$ (which is at most 2 in every case) is maximal when $G^*$ is a union of trees.

\beg{enumerate}
\item {\bf Case: $G \in \mc{S}_1^{\pr}$, and $G^* \in \mc{T}$.} The hypothesis on $G^*$ means that $G^*$ is a union of trees.
Here $\ga_G=i_G-2$. Each partition $\mu$ of $\{\la_1, \dots, \widehat{i_{G}}, \dots, \la_l\}$ into $\ga_{G}$ subsets $\mu_1, \dots, \mu_{\ga_{G}}$ contributes a term
\[
a_{\mu} C^{\mc{T}}_{\mu_1} \cdots C^{\mc{T}}_{\mu_{\ga_{G}}}.
\]
where $C^{\mc{T}}_{\mu}$ denotes the number of connected spanning trees of $K_{|\mu|}$ with indegree partition $\mu$. Here $a_{\mu}$ is the product of the following:
\beg{enumerate}
\item $\binom{d-1}{|\mu_1|,\dots,|\mu_{\ga_G}|}$
\item The number of ways of choosing $\ga_{G}$ vertices from among the connected components of $G^*$, each of which contributes 1, 2, or 3 vertices; this is 
\[
\prod_{j=1}^{\ga_G} |\mu_j|  \cdot \bigg(\sum_{1 \leq j_1<j_2 \leq \ga_G} \fr{\binom{|\mu_{j_1}|}{2}\binom{|\mu_{j_2}|}{2}}{|\mu_{j_1}| |\mu_{j_2}|} + \sum_{k=1}^{\ga_G} \fr{\binom{|\mu_k|}{3}}{|\mu_k|}\bigg).
\]
\end{enumerate}

\item {\bf Case: $G \in \mc{S}_1^{\pr}$, and $G^* \in \mc{S}_1$.} The hypothesis on $G^*$ means that exactly one connected component of $G^*$ has Betti number of 1 (i.e., a cycle), and the remaining components are trees. Here $\ga_G=i_G-1$. Each partition $\mu$ of $\{\la_1, \dots, \widehat{i_{G}}, \dots, \la_l\}$ into $\ga_{G}$ subsets $\mu_1, \dots, \mu_{\ga_{G}}$ contributes a term
\[
a_{\mu} C^{\mc{T}}_{\mu_1} \cdots C^{\mc{T}}_{\mu_{\ga_{G}}} \sum_{i=1}^{\ga_G} \fr{C^{\mc{S}_1}_{\mu_i}}{C^{\mc{T}}_{\mu_i}}.
\]
where $C^{\mc{T}}_{\mu}$ denotes the number of $|\mu|$-edged connected spanning subgraphs of $K_{|\mu|-1}$ with indegree partition $\mu$. Here $a_{\mu}$ is the product of the following:
\beg{enumerate}
\item $\binom{d-1}{|\mu_1|,\dots,|\mu_{\ga_G}|}$
\item The number of ways of choosing $\ga_{G}$ vertices, one from all but one connected component of $G^*$, which contributes 2 distinct vertices; i.e., $\prod_{j=1}^{\ga_G} |\mu_j|  \cdot \sum_{k=1}^{\ga_G} \fr{\binom{|\mu_k|}{2}}{|\mu_k|}$.
\end{enumerate}

\includegraphics[scale=.15]{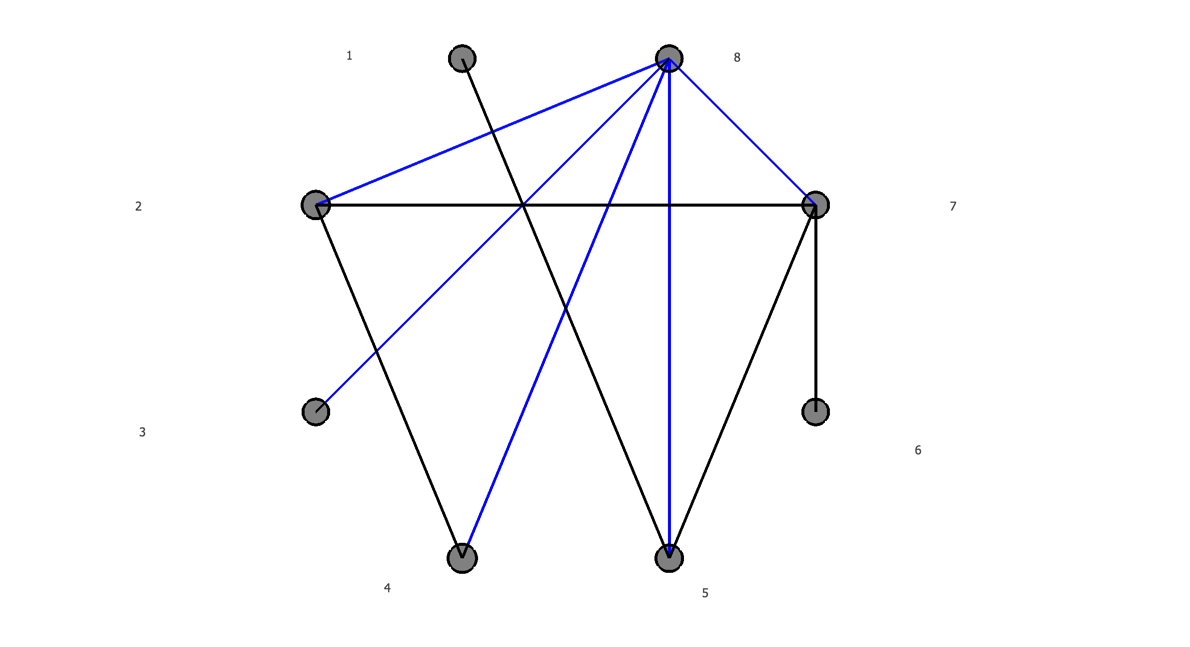}
\includegraphics[scale=.15]{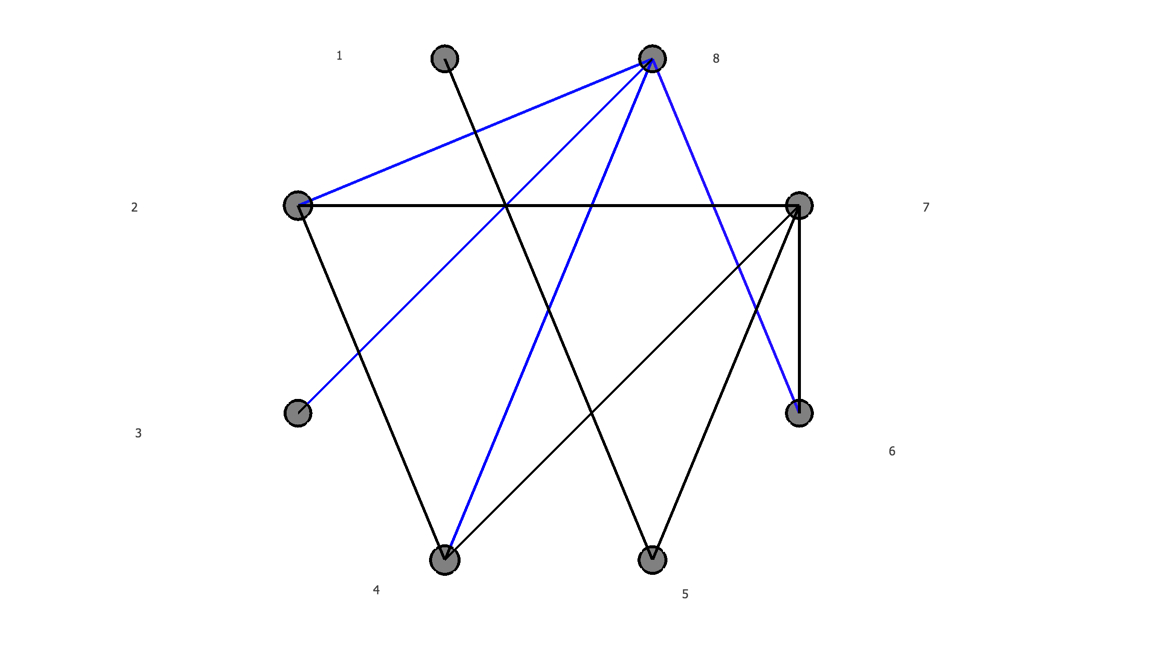}

Typical examples of Cases 1 and 2 are shown. Here $d=8$, and the edges incident to $v_d$ are highlighted in blue. Deletion of these edges yields $G^*$.

\item {\bf Case: $G \in \mc{S}_1^{\pr}$, and $G^* \in \mc{S}_1^{\pr}$.} The hypothesis on $G^*$ means that $G^*$ contains 2 linearly independent cycles, i.e., the Betti number $b_1(G^*)=2$. Here $\ga_G=i_G$, and each component of $G^*$ is linked to $v_d$ by a single edge. There are now two basic possibilities to consider: either one component of $G^*$ contains both cycles, or the cycles belong to distinct components of $G^*$. Accordingly, each partition $\mu$ of $\{\la_1, \dots, \widehat{i_{G}}, \dots, \la_l\}$ into $\ga_{G}$ subsets $\mu_1, \dots, \mu_{\ga_{G}}$ contributes terms of the form
\[
C^{\mc{T}}_{\mu_1} \cdots C^{\mc{T}}_{\mu_{\ga_{G}}} \bigg(a^{(1)}_{\mu} \sum_{i=1}^{\ga_G} \fr{C^{\mc{S}_1^{\pr}}_{\mu_i}}{C^{\mc{T}}_{\mu_i}}+ a^{(2)}_{\mu} \sum_{1 \leq j_1<j_2 \leq \ga_G} \fr{C^{\mc{S}_1}_{\mu_{j_1}}C^{\mc{S}_1}_{\mu_{j_2}}}{C^{\mc{T}}_{\mu_{j_1}}C^{\mc{T}}_{\mu_{j_2}}} \bigg).
\]
Here $a^{(1)}_{\mu}= a^{(2)}_{\mu}= \binom{d-1}{|\mu_1|,\dots,|\mu_{\ga_G}|} \cdot \prod_{j=1}^{\ga_G} |\mu_j|$.

\item {\bf Case: $G \in \mc{S}_2^{\pr}$, and $G^* \in \mc{T}$.} Here $\ga_G=i_G-1$. Each partition $\mu$ of $\{\la_1, \dots, \widehat{i_{G}}, \dots, \la_l\}$ into $\ga_{G}$ subsets $\mu_1, \dots, \mu_{\ga_{G}}$ contributes a term
\[
a_{\mu} C^{\mc{T}}_{\mu_1} \cdots C^{\mc{T}}_{\mu_{\ga_{G}}}
\]
where $a_{\mu}$ is the product of the following:
\beg{enumerate}
\item $\binom{d-1}{|\mu_1|,\dots,|\mu_{\ga_G}|}$
\item The number of ways of choosing $\ga_{G}$ vertices from among the components of $G^*$, each of which contributes 1 vertex except possibly
\beg{itemize}
\item two of these, of which one contributes a multiplicity-2 vertex, and the other contributes two multiplicity-1 vertices; or
\item a unique component  that contributes a multiplicity-2 vertex and a multiplicity-1 vertex.
\end{itemize}
The number of these possibilities is computed by
\[
2\prod_{j=1}^{\ga_G} |\mu_j|  \cdot \sum_{j=1}^{\ga_G} (|\mu_j|-1).
\]
\end{enumerate}

\item {\bf Case: $G \in \mc{S}_2^{\pr}$, and $G^* \in \mc{S}_1$.} Here $\ga_G=i_G$. Each partition $\mu$ of $\{\la_1, \dots, \widehat{i_{G}}, \dots, \la_l\}$ into $\ga_{G}$ subsets $\mu_1, \dots, \mu_{\ga_{G}}$ contributes a term
\[
a_{\mu} C^{\mc{T}}_{\mu_1} \cdots C^{\mc{T}}_{\mu_{\ga_{G}}} \sum_{i=1}^{\ga_G} \fr{C^{\mc{S}_1}_{\mu_i}}{C^{\mc{T}}_{\mu_i}}.
\]
Here $a_{\mu}$ is the product of the following:
\beg{enumerate}
\item $\binom{d-1}{|\mu_1|,\dots,|\mu_{\ga_G}|}$
\item The number of ways of choosing $\ga_{G}$ vertices, one from all but one connected component of $G^*$, which contributes a vertex with multiplicity 2; i.e., \\$\ga_G \prod_{j=1}^{\ga_G} |\mu_j|$.
\end{enumerate}

\item {\bf Case: $G \in \mc{S}_2^{\pr}$, and $G^* \in \mc{S}_2$.} The hypothesis on $G^*$ means that $G^*$ is a union of trees, exactly one of which has a multiplicity-2 edge. Components of $G^*$ are linked to $v_d$ by 1 edge each, except for one, which is linked by 2 edges. Here $\ga_G=i_G-1$.
Each partition $\mu$ of $\{\la_1, \dots, \widehat{i_{G}}, \dots, \la_l\}$ into $\ga_{G}$ subsets $\mu_1, \dots, \mu_{\ga_{G}}$ contributes a term
\[
a_{\mu} C^{\mc{T}}_{\mu_1} \cdots C^{\mc{T}}_{\mu_{\ga_{G}}} \sum_{i=1}^{\ga_G} \fr{C^{\mc{S}_2}_{\mu_i}}{C^{\mc{T}}_{\mu_i}}.
\]
Here $a_{\mu}$ is the product of the following:
\beg{enumerate}
\item $\binom{d-1}{|\mu_1|,\dots,|\mu_{\ga_G}|}$
\item The number of ways of choosing $\ga_{G}$ vertices, one from all but one connected component of $G^*$, which contributes 2 vertices; i.e., $\prod_{j=1}^{\ga_G} |\mu_j|  \cdot \sum_{k=1}^{\ga_G} \fr{\binom{|\mu_k|}{2}}{|\mu_k|}$.
\end{enumerate}

\item {\bf Case: $G \in \mc{S}_2^{\pr}$, and $G^* \in \mc{S}_2^{\pr}$.} The hypothesis on $G^*$ means that exactly one component of $G^*$ contains a cycle, and that $G^*$ contains exactly one multiplicity-2 edge. Here $\ga_G=i_G$, and each component of $G^*$ is linked to $v_d$ via a single (simple, i.e., of multiplicity 1) edge. Moreover, $G$ obeys a basic dichotomy, depending upon whether or not its multiplicity-2 edge lies along its cycle. Accordingly, each partition $\mu$ of $\{\la_1, \dots, \widehat{i_{G}}, \dots, \la_l\}$ into $\ga_{G}$ subsets $\mu_1, \dots, \mu_{\ga_{G}}$ contributes
\[
\binom{d-1}{|\mu_1|,\dots,|\mu_{\ga_G}|} \prod_{j=1}^{\ga_G} |\mu_j| \cdot \bigg(C^{\mc{T}}_{\mu_1} \cdots C^{\mc{T}}_{\mu_{\ga_{G}}} \sum_{i=1}^{\ga_G} \fr{C^{\mc{S}_2^{\pr}}_{\mu_i}}{C^{\mc{T}}_{\mu_i}}+ C^{\mc{T}}_{\mu_1} \cdots C^{\mc{T}}_{\mu_{\ga_{G}}} \sum_{1 \leq j_1 \neq j_2 \leq \ga_G} \fr{C^{\mc{S}_1}_{\mu_{j_1}}}{C^{\mc{T}}_{\mu_{j_1}}} \fr{C^{\mc{S}_2}_{\mu_{j_2}}}{C^{\mc{T}}_{\mu_{j_2}}} \bigg).
\]

\item {\bf Case: $G \in \mc{S}_3^{\pr}$, and $G^* \in \mc{T}$.} Here $\ga_G=i_G$, and all but one component of $G^*$ is linked to $v_d$ by a single simple edge; the remaining component is linked to $v_d$ by an edge of multiplicity 3. Each partition $\mu$ of $\{\la_1, \dots, \widehat{i_{G}}, \dots, \la_l\}$ into $\ga_{G}$ subsets $\mu_1, \dots, \mu_{\ga_{G}}$ contributes
\[
\binom{d-1}{|\mu_1|,\dots,|\mu_{\ga_G}|} \ga_G \prod_{j=1}^{\ga_G} |\mu_j| \cdot C^{\mc{T}}_{\mu_1} \cdots C^{\mc{T}}_{\mu_{\ga_{G}}}.
\] 
\item {\bf Case: $G \in \mc{S}_3^{\pr}$, and $G^* \in \mc{S}_3^{\pr}$.} The hypothesis on $G^*$ means that $G^*$ is a union of trees, and contains exactly one edge of multiplicity 3. Here $\ga_G=i_G$, and every component of $G^*$ is linked to $v_d$ by a unique simple edge. 
Each partition $\mu$ of $\{\la_1, \dots, \widehat{i_{G}}, \dots, \la_l\}$ into $\ga_{G}$ subsets $\mu_1, \dots, \mu_{\ga_{G}}$ contributes
\[
\binom{d-1}{|\mu_1|,\dots,|\mu_{\ga_G}|} \prod_{j=1}^{\ga_G} |\mu_j| \cdot C^{\mc{T}}_{\mu_1} \cdots C^{\mc{T}}_{\mu_{\ga_{G}}} \sum_{i=1}^{\ga_G} \fr{C^{\mc{S}_3^{\pr}}_{\mu_i}}{C^{\mc{T}}_{\mu_i}}.
\] 

\item {\bf Case: $G \in \mc{S}_{2,2}^{\pr}$, and $G^* \in \mc{T}$.} Here $\ga_G=i_G$. All but one or two components of $G^*$ are linked to $v_d$ by single simple edges; each of the remaining components is linked to $v_d$ by an edge with multiplicity 2. Each partition $\mu$ of $\{\la_1, \dots, \widehat{i_{G}}, \dots, \la_l\}$ into $\ga_{G}$ subsets $\mu_1, \dots, \mu_{\ga_{G}}$ contributes
\[
\binom{d-1}{|\mu_1|,\dots,|\mu_{\ga_G}|} \prod_{j=1}^{\ga_G} |\mu_j| \cdot C^{\mc{T}}_{\mu_1} \cdots C^{\mc{T}}_{\mu_{\ga_{G}}} \bigg(\sum_{k=1}^{\ga_G} \fr{\binom{|\mu_k|}{2}}{|\mu_k|} + \binom{\ga_G}{\ga_G-1} \bigg).
\]

\item {\bf Case: $G \in \mc{S}_{2,2}^{\pr}$, and $G^* \in \mc{S}_2$.} The hypothesis on $G^*$ means that $G^*$ is a union of trees, and contains exactly one edge of multiplicity 2. Here $\ga_G=i_G$. All but one component of $G^*$ are linked to $v_d$ by single simple edges; the remaining component is linked to $v_d$ by an edge with multiplicity 2. Each partition $\mu$ of $\{\la_1, \dots, \widehat{i_{G}}, \dots, \la_l\}$ into $\ga_{G}$ subsets $\mu_1, \dots, \mu_{\ga_{G}}$ contributes
\[
\binom{d-1}{|\mu_1|,\dots,|\mu_{\ga_G}|} \ga_G \prod_{j=1}^{\ga_G} |\mu_j| \cdot C^{\mc{T}}_{\mu_1} \cdots C^{\mc{T}}_{\mu_{\ga_{G}}} \sum_{i=1}^{\ga_G} \fr{C^{\mc{S}_2}_{\mu_i}}{C^{\mc{T}}_{\mu_i}}.
\]

\item {\bf Case: $G \in \mc{S}_{2,2}^{\pr}$, and $G^* \in \mc{S}_{2,2}^{\pr}$.} The hypothesis on $G^*$ means that $G^*$ is a union of trees, and contains two edges of multiplicity 2. Here $\ga_G=i_G$. Every component of $G^*$ is linked to $v_d$ by a single simple edge. There are two basic possibilities for $G^*$, depending upon whether or not the two multiplicity-2 edges lie along the same component. Each partition $\mu$ of $\{\la_1, \dots, \widehat{i_{G}}, \dots, \la_l\}$ into $\ga_{G}$ subsets $\mu_1, \dots, \mu_{\ga_{G}}$ contributes
\[
\binom{d-1}{|\mu_1|,\dots,|\mu_{\ga_G}|} \prod_{j=1}^{\ga_G} |\mu_j| \cdot C^{\mc{T}}_{\mu_1} \cdots C^{\mc{T}}_{\mu_{\ga_{G}}} \bigg(\sum_{i=1}^{\ga_G} \fr{C^{\mc{S}_{2,2}^{\pr}}_{\mu_i}}{C^{\mc{T}}_{\mu_i}}+ \sum_{i=1}^{\ga_G} \fr{C^{\mc{S}_{2}}_{\mu_i}}{C^{\mc{T}}_{\mu_i}} \fr{C^{\mc{S}_{2}}_{\mu_i}}{C^{\mc{T}}_{\mu_i}} \bigg).
\]

\end{enumerate}
It remains to explain how to calculate the numbers $C_{\la}^{\mc{T}}$ and $C_{\la}^{\mc{S}_i}$, $i=1, 2$. In fact, in the thesis \cite[p.29]{C2}, we conjectured that for every $d \geq 2$ and every partition $\la= (\la_1^{e_1}, \la_l^{e_l})$ of $(d-1)$,
\[
C_{\la}^{\mc{T}}= \fr{(d-1)!^2}{(d-k)! e_1! \cdots e_l! (\la_1!)^{e_1} \cdots (\la_l!)^{e_l}}
\]
where $\la_i \neq \la_j$ if $i \neq j$, and $k=\sum_{i=1}^l e_i$. This conjecture was subsequently confirmed independently (and via 2 different proofs) in \cite{DY} and \cite{SZ}.
There are no closed-form expressions known for the numbers $C_{\la}^{\mc{S}_i}$, but as usual, these may be computed recursively.

{\fl Namely}, assume that $G \in \mc{S}$ has indegree partition $\la$. Let $G^*$ denote the subgraph of $K_d$ obtained by deleting every edge of $G$ incident to $v_d$. There are four basic possibilities, depending upon whether $G$ lies in $\mc{S}_1$ or $\mc{S}_2$, 
and whether $G^*$ is a union of trees (resp., contains a doubled edge). As before, we let $\ga_G$ denote the number of connected components of $G^*$.

\beg{enumerate}
\item {\bf Case: $G \in \mc{S}_1$, and $G^* \in \mc{T}$.} Here $ \ga_G=i_G-1$. Of the components of $G^*$, all except one are linked to $v_d$ by a single edge. The remaining connected component is linked to $v_d$ along two edges.

Each partition $\mu$ of $\{\la_1, \dots, \widehat{i_{G}}, \dots, \la_l\}$ into $\ga_{G}$ subsets $\mu_1, \dots, \mu_{\ga_{G}}$ contributes a term of the form
\[
a_{\mu} C^{\mc{T}}_{\mu_1} \cdots C^{\mc{T}}_{\mu_{\ga_{G}}}
\]
to $C^{\mc{S}_1}_{\la}$. Then $a_{\mu}$ is the product of the following two quantities:
\beg{enumerate}
\item The number of ways of partitioning the set $\{1,\dots,d-1\}$ of vertices of $G^*$ into $\ga_{G}$ subsets of size $|\mu_i|, 1 \leq i \leq \ga_G$; i.e., $\binom{d-1}{|\mu_1|,\dots,|\mu_{\ga_G}|}$.

\item The number of ways of choosing $\ga_{G}$ vertices from among the connected components of $G^*$, each of which contributes 1 or 2 vertices. This is
\[
\prod_{j=1}^{\ga_G} |\mu_j|  \cdot \sum_{k=1}^{\ga_G} \fr{\binom{|\mu_k|}{2}}{|\mu_k|}.
\]
\end{enumerate}

\item {\bf Case: $G \in \mc{S}_1$, and $G^* \in \mc{S}_1$.} Here $\ga_G=i_G$. Each component of $G^*$ is linked to $v_d$ by a single edge.

Each partition $\mu$ of $\{\la_1, \dots, \widehat{i_{G}}, \dots, \la_l\}$ into $\ga_{G}$ subsets $\mu_1, \dots, \mu_{\ga_{G}}$ contributes
\[
\binom{d-1}{|\mu_1|,\dots,|\mu_{\ga_G}|} \cdot \prod_{j=1}^{\ga_G} |\mu_j| \cdot C^{\mc{T}}_{\mu_1} \cdots C^{\mc{T}}_{\mu_{\ga_{G}}} \sum_{i=1}^{\ga_G} \fr{C^{\mc{S}_1}_{\mu_i}}{C^{\mc{T}}_{\mu_i}}.
\]

\item {\bf Case: $G \in \mc{S}_2$, and $G^* \in \mc{T}$.} Here $\ga_G=i_G$. 
Each partition $\mu$ of $\{\la_1, \dots, \widehat{i_{G}}, \dots, \la_l\}$ into $\ga_{G}$ subsets $\mu_1, \dots, \mu_{\ga_{G}}$ contributes a term of the form
\[
a_{\mu} C^{\mc{T}}_{\mu_1} \cdots C^{\mc{T}}_{\mu_{\ga_{G}}}.
\]
Here $a_{\mu}$ is the product of the following two quantities:
\beg{enumerate}
\item $\binom{d-1}{|\mu_1|,\dots,|\mu_{\ga_G}|}$
\item The number of ways of choosing $\ga_{G}$ vertices, one from all but one connected component of $G^*$, which contributes a vertex with multiplicity 2; this is $\ga_G \prod_{j=1}^{\ga_G} |\mu_j|$.
\end{enumerate}

\item {\bf Case: $G \in \mc{S}_2$, and $G^* \in \mc{S}_2$.} Here $\ga_G=i_G$; each component of $G^*$ is linked to $v_d$ via a unique edge with multiplicity 1.
Accordingly, each partition $\mu$ of $\{\la_1, \dots, \widehat{i_{G}}, \dots, \la_l\}$ into $\ga_{G}$ subsets $\mu_1, \dots, \mu_{\ga_{G}}$ contributes
\[
\binom{d-1}{|\mu_1|,\dots,|\mu_{\ga_G}|} \cdot \prod_{j=1}^{\ga_G} |\mu_j|. C^{\mc{T}}_{\mu_1} \cdots C^{\mc{T}}_{\mu_{\ga_{G}}} \sum_{i=1}^{\ga_G} \fr{C^{\mc{S}_2}_{\mu_i}}{C^{\mc{T}}_{\mu_i}}.
\]

\end{enumerate}

\section{Divisor class calculations via multiple-point formulas}\label{mpt}
The secant plane divisor coefficients $P$ are not uniquely determined by the relations obtained in Section 2; rather, an additional relation is needed. In this section, we will describe an alternative method for computing secant plane divisor classes. Using this alternative approach, we are led to the following conjecture, which is borne out in every computable case.

\beg{conj}\label{conjecturalrelation}
When $r=1$, the polynomials $P_{\al}, P_{\be}, P_{\ga}, \text{ and }P_{\delta_0}$ satisfy 
\beg{equation}\label{finalrelation1}
2(d-1)P_{\al}+ (m-3)P_{\be}=(6-3g)(P_{\ga}+P_{\delta_0});
\end{equation}
when $r=s$, the polynomials $P$ satisfy
\beg{equation}\label{finalrelation2}
2(s-1)P_{\al}+ (2m-3s)P_{\be}=(6s-3m)P_{\ga}-(15m-30s+12-6g)P_{\delta_0}.
\end{equation}
\end{conj}

{\fl \bf Key observation:} Because every divisor on the stack $\mc{G}^s_m$ of curves with linear series (see \cite{Kh2} for its construction) is determined by its degrees along 1-parameter families of linear series, our general formula \eqref{basicformula} for $N_d^{d-r-1}$ determines the class of a secant plane divisor in $\mc{G}^s_m$ as a sum involving tautological classes on  $\mc{G}^s_m$. 

The upshot of the latter observation is that the conjectural relations \eqref{finalrelation1} and \eqref{finalrelation2}, coupled with the four relations among tautological coefficients $P$ obtained in Section 3, determine the classes of secant plane divisors on $\mc{G}^s_m$.

\subsection{Set-up for multiple-point formulas}
Our alternative method for calculating secant plane formulas is as follows. As before, we let $\pi: \mc{X} \ra B$ denote a one-parameter family of curves, whose total space $\mc{X}$ comes equipped with a line bundle $\mc{L}$, and whose base space $B$ comes equipped with a rank-$(s+1)$ vector bundle $\mc{V} \hra \pi_* \mc{L}$. When either $r=1$ or $r=s$ (and, conjecturally, in general), the fact that $N_{K_3}=0$ expresses the coefficient $P_{\de_0}$ in terms of the other secant plane divisor coefficients, none of which depend upon the number of singular fibers in $\pi$. Consequently, we will assume that every fiber of $\pi$ is a smooth curve. As before, the pair $(\mc{L},\mc{V})$ defines a map $f: \mc{X} \ra \mb{P} \mc{V}^*$ of $B$-schemes, whose fibers over points in $B$ are maps from curves to $s$-dimensional projective spaces. 

Now let $\mc{G}$ denote the Grassmann bundle of $(d-r-1)$-dimensional subspaces of fibers of $\mb{P} \mc{V}^*$ over $B$, and let
$\mc{I}_{\mc{X}} \sub \mc{X} \times_B \mc{G}$
denote the incidence correspondence canonically obtained from $f$. The secant plane locus of interest to us is (the pushforward to $B$ of) Kleiman's $d$th multiple-point locus \cite{Kl} associated with the projection $\rho: \mc{I}_{\mc{X}} \ra \mc{G}$. Because every fiber of $\pi$ is a smooth curve, $\pi$ is a curvilinear map. Consequently, according to \cite{Kat}, the Chow class $m_k$ of the $k$th multiple point locus of $\rho$ satisfies
\beg{equation}\label{multipleptrecursion}
m_k= \rho^* \rho_* m_{k-1}+ \sum_{i=1}^k (-1)^i p_i m_{k-i}
\end{equation}
for certain polynomials $p_i, 1 \leq i \leq k$ in the Chern classes of the virtual normal bundle $\mc{N}_{\rho}$ of $\rho$. 
Note that \cite[bot. p.11]{Kaz} gives an explicit generating series for the polynomials $p_i$.

\subsection{Evaluation of multiple-point formulas}
Evaluating the iterative formulas \eqref{multipleptrecursion} requires computing the Chern classes of the virtual normal bundle $\mc{N}_{\rho}$. These we calculate as follows. Letting $\mc{Q}_{\mc{G}}$ denote the quotient bundle on $\mc{G}$, note that because $\mc{I}_{\mc{X}} \sub \mc{X} \times_B \mc{G}$ is the zero locus of the natural map of vector bundles
\[
\mc{L}^* \ra \mc{Q}_{\mc{G}},
\]
its normal bundle $\mc{N}_{\mc{I}_{\mc{X}}/\mc{X} \times_B \mc{G}}$ is simply the pullback of $\mc{L} \otimes \mc{Q}_{\mc{G}}$ to $\mc{I}_{\mc{X}}$. On the other hand, $\mc{N}_{\mc{I}_{\mc{X}}/\mc{X} \times_B \mc{G}}$ and $\mc{N}_{\rho}$ fit into an exact sequence
\[
0 \ra \mc{T}_{\mc{X}/B} \ra \mc{N}_{\mc{I}_{\mc{X}}/\mc{X} \times_B \mc{G}} \ra \mc{N}_{\rho} \ra 0,
\]
which implies that their Chern polynomials are related by
\[
\beg{split}
c_t(\mc{N}_{\rho})&= c(-\mc{T}_{\mc{X}/B}) c(\mc{N}_{\mc{I}_{\mc{X}}/\mc{X} \times_B \mc{G}}) \\
&= (1+ t \om+ t^2 \om^2) c(\mc{N}_{\mc{I}_{\mc{X}}/\mc{X} \times_B \mc{G}})
\end{split}
\]
where $\om= c_1(\om_{\mc{X}/B})$ is the first Chern class of the relative dualizing sheaf of $\pi$.

It is now a relatively straightforward matter to write a computer program to compute secant plane divisor classes; this we have done in Maple \cite{Ma}. Code is provided on the webpage \url{www.mast.queensu.ca/~cotteril}.

\subsection{Examples}\label{gsmexamples}
In this subsection, we record several secant plane formulas in cases where  either $r=1$ or $r=s$. Note that with the exception of the first two and the fifth, which Ran \cite{R2} has computed using his intersection theory for families of rational curves, these are new. Formulas in the range $d>6$ occupy too much space to included here.

\beg{itemize}
\item {\bf $r=1, d=2, s=3$ (case of 3-dimensional series with double points).} Here
{\small
\[
2! N_2^0= (-6 + 2 m)\al - 4\be + (2 g - 2 + 3 m - m^2) c - \ga + \de_0.
\]
}
\item {\bf $r=1, d=3, s=5$ (case of 5-dimensional series with trisecant lines).} Here
{\small
\[
\beg{split}
3! N_3^1&= (3m^2-27m-6g+66)\al+(72-12m)\be+(28-3m)\ga+ (3m-20) \de_0 \\
&+(24-m^3+9m^2+6mg-26m-24g)c.
\end{split}
\]
}
\item {\bf $r=1, d=4, s=7$ (case of 7-dimensional series with 4-secant 2-planes).} Here
{\small
\[
\beg{split}
4! N_4^2&= (-1008 + 168 g - 24 m g - 72 m^2  + 452 m + 4 m^3) \al+ (360 m - 1440 + 48 g - 24 m^2) \be+ \\
&+ (372 g -360 + 342 m - 119 m^2  - m^4  + 18 m^3  - 12 g^2  - 132 m g + 12 m^2  g) c \\
&+(12 g -720 + 130 m - 6 m^2) \ga +(6 m^2 - 98 m - 12 g + 432) \de_0. 
\end{split}
\]
}
\item {\bf $r=1, d=5, s=9$ (case of 9-dimensional series with 5-secant 3-planes).} Here
{\small
\[
\beg{split}
5! N_5^3&= (1020mg - 60m^2g - 4500g + 60g^2  + 19560 + 5m^4  + 1735m^2  - 150m^3  - 9270m) \al\\
&+ (240mg - 2400g + 33600 - 40m^3  - 10160m + 1080m^2) \be \\
&+ (20000 + 60mg - 800g + 370m^2  - 10m^3  - 4640m) \ga \\
&+(20 m^3  g - 60 mg^2- 420 m^2  g + 6720 + 480 g^2  + 2980 m g - 5944 m + 30 m^4 \\
&- 355 m^3  + 2070 m^2  - m^5  - 7200 g) c \\
&+(60mg + 640g + 10m^3  + 2960m - 290m^2  - 10720) \de_0.
\end{split}
\]
}
\item {\bf $r=2, d=3, s=2$ (case of 2-dimensional series with triple points).} Here
{\small
\[
\beg{split}
3! N_3^0&=(3m^2-18m-6g+30) \al+ (18-3m) \be+ 4 \ga -2 \de_0 \\
&+ (12m^2-2m^3+6mg- 22m+12 -12g) c.
\end{split}
\]
}
\item {\bf $r=3, d=5, s=3$ (case of 3-dimensional series with 5-secant lines).} Here
{\small
\[
\beg{split}
4!N_5^1&=(10m^4  - 180m^3  + 1250m^2  + 5160 - 60m^2g - 4020m + 600mg + 60g^2  - 1620g) \al\\
&+ (360m^2  - 20m^3  + 4800 + 60mg - 2200m - 480g) \be+ (1520 - 450m + 40m^2  - 80g) \ga \\
&+ (2400 + 2190 m^2  + 1940mg - 2640 g - 635 m^3  - 60 m g^2  - 480 m^2  g \\
&+ 40 m^3  g - 3680 m - 5 m^5+ 90 m^4  + 240 g^2) c + (40g-20m^2 + 210m - 640) \de_0.
\end{split}
\]
}
\end{itemize}
Note that in each of the above examples, we have realized the class of a divisor on the stack $\mc{G}^s_m$ in terms of tautological classes $c$, $\al$, $\be$, $\ga$, and $\de_0$. 

\section{Le Barz's cycle-theoretic secant planes}\label{lebarz}
In \cite{Lb1} and \cite{Lb2}, P. Le Barz adopts a different approach to the enumeration of secant planes of (fixed) projective curves, based on the theory of excess intersection. Here we review his method, recasting it in slightly more generality, compare it to our own, and suggest how one might adapt it to the setting of one-parameter families of curves with linear series.

\subsection{Cycles in the Grassmannian associated to secant planes}
Let $n \geq 1$, and assume that $C \sub \mb{P}^n$ is a smooth curve of degree $m$ and genus $g$. Fix $d^{\pr}<n$, and set $\mb{G}:= \mb{G}(d^{\pr},n)$. Requiring a $d^{\pr}$-plane to intersect $C$ is a codimension-$(n-d^{\pr}-1)$- condition; accordingly, those $d^{\pr}$-planes that have degree-$k$ intersections with $C$ have expected codimension $k(n-d^{\pr}-1)$ in $\mb{G}$. While it is not necessarily the case that the actual and expected codimensions agree, nevertheless there is always a well-defined $k$-secant {\it cycle} class $\mbox{Sec}_k(C) \in A^{k(n-d^{\pr}-1)}(\mb{G})$, obtained as follows.

{\fl Let} $\mc{S}$ denote the tautological rank-$(d^{\pr}+1)$ subbundle of $\mb{G}$, and let $\mc{T}= \mb{P}(\mc{S}) \hra \mb{G} \times \mb{P}^n$ be the corresponding projective bundle. Then
\[
\mbox{Co}^k(d^{\pr},n):= \mc{T}^{[k]}
\]
is the incidence scheme of pairs $(\La, Y)$ where $Y \sub \La$ is a scheme of $k$ coplanar points.
The incidence scheme comes equipped with projections $s$ and $t$ to $\mb{G}$ and $(\mb{P}^n)^{[k]}$, respectively. Here $s: \mbox{Co}^k(d^{\pr},n) \ra \mb{G}$ is a fibration with $d^{\pr}k$-dimensional fiber $(\mb{P}^{d^{\pr}})^{[k]}$. 
Whence,
\[
\dim \mbox{Co}^k(d^{\pr},n)= d^{\pr}k+ (d^{\pr}+1)(n-d^{\pr}).
\]
The incidence scheme contains a distinguished subscheme
\[
Z:= t^{-1}(C^{[k]})
\]
made up of pairs $(\La,Y)$ with $Y \sub \La \cap C$. Note that the embedding of $C$ in $\mb{P}^n$ canonically induces a codimension-$(nk-k)$ embedding of $C^{[k]}$ in $(\mb{P}^n)^{[k]}$, for every positive integer $k \geq 1$. The expected dimension of $Z$ is thus
\[
\beg{split}
\mbox{exp. dim }Z&=\dim \mbox{Co}^k(d^{\pr},n)- (nk-k) \\
&=(d^{\pr}+1)(n-d^{\pr})+ k-k(n-d^{\pr}).
\end{split}
\]
On the other hand, when $k$ is sufficiently large relative to $m$ and $n$, we have
\[
Z \cong C^{[k]} \times \mb{G}(C)
\]
where $\mb{G}(C) \sub \mb{G}$ is the subvariety of $d^{\pr}$-planes containing $C$. In particular, the difference between the actual and expected dimensions of $Z$ is
\beg{equation}\label{excessdimensions}
\beg{split}
e&=k+ (d^{\pr}-\ga)(n-d^{\pr})- [(d^{\pr}+1)(n-d^{\pr})+k-k(n-d^{\pr})] \\
&=(k-\ga-1)(n-d^{\pr})
\end{split}
\end{equation}
where $\ga= \dim \Ga$ is the dimension of the linear span $\Ga$ of $C$. Note that $\mb{G}(C)=\mb{G}(\Ga)$.

{\fl Finally, let $j$ denote the inclusion of $Z$ in $\mbox{Co}^k(d^{\pr},n)$, and set}
\beg{equation}\label{secantcycle}
\mbox{Sec}^{d^{\pr},n}_k(C):= s_* j_*(C^{[k]} \cdot \mbox{Co}^k(d^{\pr},n)).
\end{equation}

\subsection{The case $r=1$}
We have
\[
Z_{g,m}(z)= \sum_{d \geq 0} \deg(\mbox{Sec}^{d-2,2d-2}_d(C)) t^d.
\]

Le Barz's determination of $Z_{g,m}(z)$ is based on 2 principles, namely:
\beg{enumerate}
\item $Z_{g,m}(z)$ satisfies a multiplicativity property; namely,
\beg{equation}\label{multiplicativity}
Z_{g_1+g_2-1,m_1+m_2}(z)= Z_{g_1,m_1}(z) \cdot Z_{g_2,m_2}(z).
\end{equation}
\item For any {\it particular} values of $g$ and $m$, $Z_{g,m}(z)$ may be computed using Fulton's excess intersection formula, in tandem with Schubert calculus.
\end{enumerate}

{\fl \bf Multiplicativity.} Let $C_1$ and $C_2$ be disjoint smooth curves in $\mb{P}^{2d-2}$ of degrees $m_1, m_2$ and genera $g_1, g_2$, respectively. By cleverly specializing the relative positions of $C_1$ and $C_2$, Le Barz shows \cite[Lemme 2]{Lb2} that the number of $d$-secant $(d-2)$-planes to $C_1 \cup C_2$ is given by
\[
N_d= \sum_{i=0}^d N_i N_{d-i},
\]
which immediately implies \eqref{multiplicativity}.

{\fl \bf Determination of $Z_{g,m}(z)$ for particular values of $g$ and $m$.} Let $C \st{f}{\hra} \mb{P}^n$ be a smooth curve of genus $g$, embedded via a linear series of degree $m$. Correspondingly, for every integer $k \geq 1$, there is an induced inclusions of Hilbert schemes $f^{[k]}: C^{[k]} \hra (\mb{P}^n)^{[k]}$.  As explained in \cite{Lb1}, for $k >>0$, the excess intersection formula realizes the intersection product \eqref{secantcycle} as the $e$th Chern class of the virtual normal bundle
\[
\mc{F}:= t^* \mc{N}_{f^{[k]}}- \mc{N}_{Z/\mbox{Co}^k(d,n)}
\]
where $e$ is defined as in \eqref{excessdimensions}.

{\fl \bf The fundamental exact sequence}. Let $\mc{T}_k \sub C^{[k]} \times C$ denote the tautological incidence correspondence, equipped with its natural projections $p_k$ and $q$ to $C^{[k]}$ and $C$, respectively. Moreover, let
\[
\phi_{k+1}: C^{[k]} \times C \ra C^{[k+1]}
\]
denote the map defined by summation of cycles.

A result of Ran \cite{R1} establishes that the bundles $\mc{N}_{f^{[k]}}$ relate to one another via the following short exact sequences, for all $k \geq 1$:
\beg{equation}\label{normalbdlseq}
0 \ra q^* \mc{N}_f \ot \mc{O}(-\mc{T}_k) \ra \phi_{k+1}^* \mc{N}_{f^{[k+1]}} \ra p_k^* \mc{N}_{f^{[k]}} \ra 0
\end{equation}

Note that $\mc{N}_f=f^* \mc{T}_{\mb{P}^n}/\mc{T}_C$. Whence, an easy application of the Euler sequence for $\mb{P}^n$ yields
\[
c_t(\mc{N}_f)= 1+ \mc{L} t
\]
where $\mc{L}=[(n+1)m+ 2g-2] \{\text{pt}\}$.

{\fl \bf The Chern polynomial of $\mc{N}_{f^{[k]}}$ splits linearly.}
It is natural to pull back the fundamental exact sequence \eqref{normalbdlseq} to the Cartesian product $C^{k+1}$. In doing so, we adopt the following notational convention: for all $l \geq k$, the symmetric product $C^{[k]}$ in \eqref{normalbdlseq} is the quotient of the first $k$ copies of $C$ in the Cartesian product $C^l$. Doing so allows us to unambiguously omit pullbacks. The sequence \eqref{normalbdlseq} pulls back to the following sequence on $C^{k+1}$:
\beg{equation}\label{pbseq}
0 \ra  \mc{N}_f \ot \mc{O}(-\De_{k+1}) \ra \mc{N}_{f^{[k+1]}} \ra \mc{N}_{f^{[k]}} \ra 0
\end{equation}
where $\De_{k+1}$ is the ``large" diagonal comprising $(k+1)$-tuples $(x_i)_{i=1}^{k+1}$ such that $x_j=x_{k+1}$ for some index $1 \leq j \leq k$. It follows immediately that on $C^k$, the Chern polynomial of $\mc{N}_{f^{[k]}}$ splits linearly in the same fashion as the $k$th secant bundle:
\beg{equation}\label{Nfksplit}
c_t(\mc{N}_{f^{[k]}})= \prod_{j=1}^k [1+ (\mc{L}-\De_j)t]
\end{equation}
where $\De_1:=0$. Indeed, it is a conceivable that the graph-theoretic interpretation of $N_d$, and of tautological coefficients (in the setting of one-parameter families of curves) described in Subsection~\ref{graphtheoretic} can be explained on the basis of Le Barz's cycle-theoretic approach, though we haven't pursued this.

{\fl \bf Preliminaries regarding the Chern classes of $\mc{N}_{Z/\mbox{Co}^k(d,n)}$}.
Let $\Ga$ denote the linear span of $C$, and let $\mb{G}(C)=\mb{G}(\Ga) \sub \mb{G}(d,n)$ denote the subvariety of $d$-planes containing $C$. Let $W:=\mbox{Co}^k(d,n)$, and $Y:=s^{-1}\mb{G}(C)$. Note that $Z$ is contained in $Y$, and we have
\[
c_t(\mc{N}_{Z/W})= c_t(\mc{N}_{Z/Y}) c_t(\mc{N}_{Y/W})|_Z.
\]

{\fl Moreover}, because $s$ is a fibration,
\[
\mc{N}_{Y/W}= s^{\pr,*} \mc{N}_{\mb{G}(X)/\mb{G}(d,n)} = s^{\pr,*} \mc{N}_{\mb{G}(\Ga)/\mb{G}(d,n)}
\]
where $s^{\pr}= s|_Y$.

{\fl Now} note that
\beg{equation}\label{Nga}
\mc{N}_{\mb{G}(\Ga)/\mb{G}(d,n)}= (\mc{Q}^{\pr})^{\oplus \ga+1}
\end{equation}
where $\mc{Q}^{\pr}$ is the tautological quotient bundle on $\mb{G}(C)$, and $\ga= \dim(\Ga)$. Indeed, when $\Ga$ is a point, \eqref{Nga} follows easily from the standard description of the tangent space to $\mb{G}(d,n)$, and the exact sequence that realizes $\mc{N}_{\mb{G}(\Ga)/\mb{G}(d,n)}$ as a quotient of $\mc{T}_{\mb{G}(d,n)}$. In general, \eqref{Nga} follows inductively, using a filtration of $\Ga$ by a complete flag of linear subspaces.

{\fl It} follows that
\beg{equation}\label{NYW}
c_t(\mc{N}_{Y/W})= s^{\pr,*} c_t(\mc{Q}^{\pr})^{\oplus \ga+1}.
\end{equation}

{\fl The} Chern polynomial of $\mc{N}_{Z/Y}$ may be computed as follows. 
Let $\wt{f}: \mb{G}(\Ga) \times X \ra \mb{P}(\mc{S}|_{\mb{G}(\Ga)})$ denote the canonical inclusion. Whenever $k$ is large, we have $\mc{N}_{Z/Y}= \mc{N}_{\wt{f}^{[k]}}$, where
\[
\wt{f}^{[k]}: \mb{G}(\Ga) \times C^{[k]} \ra \mb{P}(\mc{S}^*|_{\mb{G}(\Ga)})^{[k]}
\]
is the canonically-induced inclusion of Hilbert schemes. The obvious generalizations of \eqref{normalbdlseq} and \eqref{pbseq} allow $c_t(\mc{N}_{\wt{f}^{[k]}})$ to be computed inductively, given $c_t(\mc{N}_{\wt{f}})$. 

 {\fl \bf An auxiliary exact sequence, and a key technical lemma for normal bundles.} The calculation of $\mc{N}_{f^{[k]}}$ and $\mc{N}_{\wt{f}^{[k]}}$ fits into the following more general framework. Namely, let $B$ be a variety, and let $E$ be a rank-$(n+1)$ vector bundle over $B$. Let $\Sig=\mb{P}(E^*)$ be the corresponding $\mb{P}^n$-bundle. Let $\wt{\Ga}$ be a fixed $(\ga+1)$-dimensional vector space, with $\Ga=\mb{P}(\wt{\Ga}^*)$. Here $B \times \Ga \sub \Sig$ is a trivial subbundle. Now let $C \sub \Ga$ denote a fixed smooth curve of degree $m$ and genus $g$ with linear span $\Ga$. We have embeddings
\[
\wt{f}^{[k]}: B \times C^{[k]} \hra (\Sig/B)^{[k]}
\]
for all $k \geq 1$, where $(\Sig/B)^{[k]}$ denotes the $k$th relative Hilbert scheme of the $B$-scheme $\Sig$.

{\fl In} \cite{Lb1}, Le Barz computes $\mc{N}_{\wt{g}^{[k]}}$ in the cases where $C$ is either a line or a conic; because of the multiplicativity property \eqref{multiplicativity} of $Z_{g,m}(z)$, these two cases suffice. To do so, he uses an immediate generalization of the fundamental exact sequence \eqref{normalbdlseq} that relates $\mc{N}_{\wt{f}^{[k]}}$ and $\mc{N}_{\wt{f}^{[k+1]}}$, namely:
\beg{equation}\label{normalbdlseqbis}
0 \ra q^* \mc{N}_{\wt{f}} \ot \mc{O}(-\mc{T}_k) \ra \phi_{k+1}^* \mc{N}_{\wt{f}^{[k+1]}} \ra p_k^* \mc{N}_{\wt{f}^{[k]}} \ra 0
\end{equation}
He also makes use of two crucial additional inputs. The first of these is the fact that the classes of the tautological divisor and the pullback via the summation map $\phi^{[k+1]}: C^{[k]} \times C \ra C^{[k+1]}$ of the hyperplane class $H_k$ on $C^{[k]} \cong \mb{P}^k$ are given, respectively, by
\beg{equation}\label{1stinput}
\mc{T}_k= P_k+ kQ, \text{ and } \phi_k^* H_k= P_k+Q
\end{equation}
where $P_k= p_k^* H_k$, and $Q= q^* \{\mbox{pt}_C\}$.

The second input is the following short exact sequence on $C \cong \mb{P}^1$, valid for all $k \geq 1$:
\beg{equation}\label{auxses}
0 \ra \mc{O}(-(k-1)) \ra \mc{O}^{\oplus k} \ra \mc{O}(1)^{\oplus k-1} \ra 0.
\end{equation}

{\fl The} auxiliary sequence \eqref{auxses} becomes useful when pulled back to $C^{[k]} \times C \times B$, where it implies that for all $k \geq 1$ and for every vector bundle $\Xi$ on $C^{[k]} \times C \times B$,
\beg{equation}\label{auxses2}
0 \ra \Xi(-(k-1)Q) \ra \Xi^{\oplus k} \ra \Xi(Q)^{\oplus k-1} \ra 0
\end{equation}
is exact.
Applying \eqref{auxses2}, \eqref{normalbdlseq}, and the first input \eqref{1stinput} exactly as in the proof of \cite[Prop. 3]{Lb1}, we deduce the following result for degree-$m$ rational curves.
\beg{thm}
For all $k \geq 1$, $m>1$,
\beg{equation}\label{arbdegrational}
c_t(\mc{N}_{\wt{f}^{[k]}})= c_t(\hat{E})^{m+1} \cdot c_t(\hat{E} \ot \mc{O}(-H_k))^{k-m-1} \cdot (1-H_k t)^{k+1-(\ga+1)m}
\end{equation}
where $\hat{E}= E/(B \times \wt{\Ga})$. 
\end{thm}

The formula \eqref{arbdegrational} generalizes Le Barz's \cite[Prop. 3bis]{Lb1}. In fact, however, we can go further. Namely, assume that $g>0$. In that case, the analogue of the hyperplane class $H_k$ is played by the class $x_k$ of the locus $X_{k, \put(1,0){\framebox(2,2)}} \hspace{2pt} \sub C^{[k]}$ comprising $k$-cycles whose support contains a fixed point in $C$. Indeed, as explained in \cite[Ch. 7, Prop. 2.1]{ACGH}, $\mc{O}(X_{k, \put(1,0){\framebox(2,2)}}\hspace{2pt})$ is canonically isomorphic to the polarization $\mc{O}(1)$ on the projective bundle given by $C^{[k]} \ra \mbox{Pic}^k(C)$ whenever $k \geq 2g-1$.

Accordingly, applying \cite[p. 338]{ACGH}, we find that for $g \geq 1$, the analogue of the first input \eqref{1stinput} is given by
\beg{equation}\label{1stinputbis}
\mc{T}_k= P_k+ \ga+ kQ, \text{ and } \phi_k^* H_k= P_k+Q
\end{equation}
where $P_k= p_k^* x_k$, $Q$ is as before, and $\ga$ is the $(1,1)$-part of the class of the diagonal $\De \sub C \times C$. Recall that the $(1,1)$-parts of $H^*(C^{[k]} \times C, \mb{Q})$ and $H^*(C \times C, \mb{Q})$ are isomorphic for all $k \geq 1$.

{\fl Similarly}, it is natural to ask for analogues of the auxiliary short exact sequences \eqref{auxses} and \eqref{auxses2}. For the purposes of intersection theory, however, we need only the resulting Chern polynomial identity
\beg{equation}\label{chernpolyid}
c_t(\Xi)^k= c_t(\Xi(-(k-1)Q) c_t(\Xi(Q))^{k-1},
\end{equation}
which holds irrespective of the genus of $C$.

{\fl It} now follows by the same argument used in the case of rational curves that
\[
c_t(\mc{N}_{\wt{f}^{[k]}})= c_{\ga} \cdot c_t(\hat{E})^{m+1} \cdot c_t(\hat{E} \ot \mc{O}(-H_k))^{k-m-1} \cdot (1-H_k t)^{k+1-2g-(\ga+1)m}.
\]
where $c_{\ga}$ is a contribution arising from $\ga$.
For $k=1$, $c_{\ga}=1$. When $k=2$, a straightforward calculation using \eqref{normalbdlseqbis}, \eqref{1stinputbis}, \eqref{chernpolyid}, and the splitting principle yields $c_{\ga}=(1-\ga t)^{n+1}$. More generally, we have $c_{\ga}= 
(1-\ga t)^{(k-1)(n+1)}$, from which we deduce the following result, valid for embedded curves of positive genus.
\beg{thm}
For all $k \geq 1$,
\beg{equation}\label{posgenusfmla}
c_t(\mc{N}_{\wt{f}^{[k]}})= (1-\ga t)^{(k-1)(n+1)} \cdot c_t(\hat{E})^{m+1} \cdot c_t(\hat{E} \ot \mc{O}(-H_k))^{k-m-1} \cdot (1-H_k t)^{k+1-2g-(\ga+1)m}.
\end{equation}
\end{thm}

{\bf \fl Extending the cycle-theoretic method to one-parameter families of curves.}
A nice project would be to generalize the results described above to the setting of one-parameter families of curves. A potential application would be a proof of Conjecture 3, by applying a generalization of \eqref{posgenusfmla} and its consequences to a suitably chosen collection (varying with the incidence parameter $d$) of non-isotrivial families of smooth curves. Note that the fibers of such a family will necessarily be of positive genus; otherwise the Hodge class $\ga$ will evaluate to zero along along the family. On the other hand, we may suppose that $s$, the dimension of the ambient projective space, remains fixed as $d$ increases; this is a key feature of the cycle-theoretic method. A final point is that when $r=1$, the classes of the Schubert monomials that appear in the calculation of $c_e(\mc{F})$ are each point classes, i.e., of degree 1 \cite[Lemme 1]{Lb2}. (The Chern classes of the tautological quotient bundle on the Grassmannian are special Schubert cycles, which explains their appearance here.) This is an important particularity of the $r=1$ case, which we expect continues to hold in the relative setting.

\section{Secant plane divisors on $\mc{G}^s_m$ and on $\ov{\mc{M}}_g$}
Thus far, we have seen how to determine the coefficients $P_{\al}, P_{\be}, P_{\ga}, \text{ and } P_{\de_0}$ of secant plane divisors on the space of linear series $\mc{G}^s_m$. For the sake of calculation, we have assumed $\rho=0$; whenever this is the case, every secant plane divisor pushes forward to a divisor on $\ov{\mc{M}}_g$. Khosla's determination of the Gysin map \cite{Kh2}, which we review now, will allow us to compute the coefficients of the Hodge class $\la$ and of the ``irreducible" boundary divisor $\de_0$, of secant plane divisors on $\ov{\mc{M}}_g$. 

\subsection{Recapitulation of Khosla's work}\label{dprecap}
Let $\widetilde{\mc{M}}_{g,1}$ denote the open substack of $\ov{\mc{M}}_{g,1}$ equal to the complement of the closure of the substack swept out by reducible unions of smooth curves intersecting transversely in two points. Let $\widetilde{\pi}: \mc{C} \ra \widetilde{\mc{M}}_{g,1}$ denote the universal curve, with relative dualizing sheaf $\widetilde{\om}$. Recall that for all $g \geq 3$ \cite{Ha},
\[
\mbox{Pic}(\widetilde{\mc{M}}_{g,1}) \ot \mb{Q}= \mb{Q} \la \op \mb{Q} \de_0 \op_{i=1}^{g-1} \mb{Q} \de_i \op \mb{Q} \psi.
\]
Here 
\[
\la= c_1(\widetilde{\pi}_* \widetilde{\om}) \text{ and } \psi=c_1(\om_{\widetilde{\mc{M}}_{g,1}/\widetilde{\mc{M}}_g}),
\]
while $\de_0$ corresponds to irreducible nodal curves, and $\de_i, i \geq 1$ corresponds to reducible unions of curves of genera $i$ and $(g-i)$ marked along the component of genus $i$.

There is, correspondingly, a Deligne-Mumford stack $\mc{G}^s_m$ of curves with linear series. When $\rho$ is nonnegative, a unique component of the stack of linear series, which we denote abusively by $\mc{G}^s_m$, dominates the moduli stack. Moreover, the projection $\eta: \mc{G}^s_m \ra \widetilde{\mc{M}}_{g,1}$ is generically smooth with fiber dimension $\rho$.

Now let $\pi: \mc{C}^s_m \ra \mc{G}^s_m$ denote the universal curve. There is a coherent sheaf $\mc{L}$ on $\mc{C}^s_m$ with torsion-free fibers, whose degree is $m$ on the marked component of every fiber, and whose degree is zero on unmarked components of fibers. Furthermore, $\mc{L}$ is trivialized along the marked section of $\pi$. It is not hard to see that the preceding two properties characterize $\mc{L}$ uniquely. Finally, there is a subbundle
\[
\mc{V} \ra \pi_* \mc{L}
\]
whose fibers are marked aspects of linear series.

In \cite[Thm. 2.11]{Kh2}, Khosla computes the images under the Gysin pushforward
\[
\eta_*: A^1(\mc{G}^s_m) \ra A^1(\widetilde{\mc{M}}_{g,1})
\]
of the tautological classes
\[
\al= \pi_*(c_1^2(\mc{L})), \be= \pi_*(c_1(\mc{L}) \cdot c_1(\om)), \text{ and } c=c_1(\mc{V})
\]
where $\om=\om_{\mc{C}^s_m/\mc{G}^s_m}$ is the relative dualizing sheaf. Note that $\al$, $\be$, and $c$ are precisely those tautological classes that appear in the basic secant plane formula \eqref{basicformula}. Moreover, there is no mention of the standard class $\ga= \pi_*(c_1^2(\om))$ here; that is because, as explained in \cite[(3.110)]{HM},
\beg{equation}\label{gammarelation}
\ga= 12 \la -\de_0.
\end{equation}

For our purposes, the contributions of $\psi$ and of $\de_i, i \geq 1$ to the pushforwards of the standard classes are immaterial, so we omit them. Khosla's formulas, streamlined in this way, read as follows.
\beg{equation}\label{Gysin}
\beg{split}
\eta_* \al&= mN\biggl[\fr{(gm-2g^2+8m-8g+4)}{(g-1)(g-2)} \la+ \fr{2g^2-gm+3g-4m-2}{6(g-1)(g-2)} \de_0\biggr], \\
\eta_* \be&= mN \biggl[ \fr{6}{g-1} \la- \fr{1}{2(g-1)} \de_0 \biggr], \text{ and } \\
\eta_* c&= N \biggl[ \fr{-(g+3) \xi+5s(s+2)}{2(g-1)(g-2)} \la+ \fr{(g+1) \xi- 3s(s+2)}{12(g-1)(g-2)} \de_0 \biggr]
\end{split}
\end{equation}
where
\beg{equation}\label{xi,n}
\beg{split}
N &= \fr{g! \cdot \prod_{i=1}^s i!}{\prod_{i=0}^s (g-m+s+i)!} \text{ is the degree of the covering }\eta, \text{ and }\\
\xi &= 3(g-1)+ \fr{(s-1)(g+s+1)(3g-2m+s-3)}{g-m+2s+1}.
\end{split}
\end{equation}

Using the equations \eqref{explicitP}, \eqref{explicitPbis}, \eqref{gammarelation}, and \eqref{Gysin}, we can explicitly determine the class $\mbox{Sec}$ of any secant plane divisor on $\ov{\mc{M}}_g$, modulo the boundary classes $\de_i, i \geq 1$ whenever $r=1$ or $r=s$. 
Namely, we have, modulo contributions from boundary divisors corresponding to reducible curves,
\beg{equation}\label{SecMg}
\beg{split}
\mbox{Sec}&= P_{\al} \eta_* \al+ P_{\be} \eta_* \be+ P_c \eta_* c+ P_{\ga} \cdot N(12 \la-\de_0) + N P_{\de_0} \de_0 \\
&=b_{\la} \la- b_0 \de_0
\end{split}
\end{equation}
where $b_{\la}=b_{\la}(d)$ and $b_0=b_0(d)$ are explicitly determined rational functions of $g$ and $m$, for any given choice of $d$. 

\subsection{Slope calculations}\label{slopexamples}
Recall \cite{HM2} that the slope of an effective divisor $D \sub \ov{\mc{M}}_g$ with class 
\beg{equation}\label{stdexpansion}
D= b_{\la} \la - b_0 \de_0 - \sum_{i=1}^{\lfl \fr{g}{2} \rfl} b_i
\end{equation}
is defined to be the quantity
\[
\mbox{slope}(D)= \fr{b_{\la}}{\min_i \{b_i\}}.
\]
As explained in \cite[Cor. 1.2]{FP}, we have $\mbox{slope}(D)= \fr{b_{\la}}{b_0}$ whenever $g \leq 23$ and provided
\beg{equation}\label{threshold}
\beg{split}
&\fr{b_{\la}}{b_0} \leq 6+ \fr{11}{\lfl \fr{g}{2} \rfl+1}, \text{ and}\\
&\fr{b_{\la}}{b_0} \leq \fr{88828}{12870} \text{ whenever }20 \leq g.
\end{split}
\end{equation}
We have checked that the ratio $\fr{b_{\la}}{b_0}$ of the first two coefficients in the expansion \eqref{stdexpansion} of $\mbox{Sec}$ in terms of standard classes satisfies the conditions \ref{threshold} whenever $g \leq 23$ and either $r=1$ or $r=s$. It follows that whenever $r=1$ or $r=s$,
\beg{equation}\label{slope}
\mbox{slope}(\mbox{Sec})= \fr{b_{\la}}{b_0}
\end{equation}
for all $g \leq 23$.
We expect, moreover, that the equation \eqref{slope} holds for {\it all} $g$. In the following table, we compile slopes of some secant plane divisors in the case $r=1$.

\medskip
\beg{tabular}{| r | r | r | r | r | r | r | r |} \hline
Genus $g$ & $d$ & $s$ & $m$ & $\fr{b_{\la}}{b_0}- (6+\fr{12}{g+1})$ & $\fr{b_{\la}}{b_0}- (6+ \fr{11}{\lfl \fr{g}{2} \rfl+1})$ & $\fr{b_{\la}}{b_0}- \fr{8828}{12870}$\\
\hline \hline
8 & 2 & 3 & 9 & 0 & $-13/15$ & N/A\\
\hline
12 & 2 & 3 & 12 & $693/12389$ & $-3952/6671$ & N/A\\
\hline
16 & 2 & 3 & 15 & $756/13379$ & $-3257/7083$ & N/A\\
\hline
20 & 2 & 3 & 18 & $1539/30247$ & $-1632/4321$ & $-7775369/27805635$\\
\hline
12 & 3 & 5 & 15 & $308/6539$ & $-2117/3521$ & N/A\\
\hline
18 & 3 & 5 & 20 & $32232/596239$ & $-130031/313810$ & N/A\\
\hline
16 & 4 & 7 & 16 & $2520/46427$ & $-11357/24579$ & N/A\\
\hline
20 & 5 & 9 & 20 & $2508/47159$ & $-2529/6737$ & $-12023068/43352595$ \\
\hline
\end{tabular}

\medskip
Note that all entries in the second-to-last column are negative, as are all entries in the last column in all cases where $g \geq 20$. It follows that $\fr{b_{\la}}{b_0}$ computes the slope in every case listed. On the other hand, the fact that
\[
\fr{b_{\la}}{b_0}- \biggl(6+\fr{12}{g+1}\biggr) \geq 0
\]
in every case shows that in each case, the slope of $\mbox{Sec}$ is at least that of the Brill--Noether divisor on $\ov{\mc{M}}_g$. The lone zero at the top of column 5 is explained by the fact that every curve of genus 8 that admits a $g^3_9$ with nodes also carries a $g^2_7$, and $\rho(8,2,7)=-1$. So the corresponding secant plane divisor is a Brill--Noether divisor on $\ov{\mc{M}}_8$.

\subsection{Nonemptiness of secant plane divisors with $r=1$}\label{nonemptyr=1}
In this section, we prove the following result.
\beg{thm}\label{nonemptysecant}
Secant plane divisors on $\ov{\mc{M}}_g$ are nonempty whenever $\rho=0$ and $r=1$.
\end{thm}

\beg{proof}
Since the pushforward $A^1(\mc{G}^s_m) \ra A^1(\ov{\mc{M}}_g)$ is finite and nonzero, it suffices to show that the classes of the corresponding secant plane divisors on $\mc{G}^s_m$ are nonzero. Moreover, in light of the calculation carried out for the first test family in Section~\ref{families}, the desired nonvanishing property will follow from showing that the tautological coefficient $P_c=P_c(d,g,m)$ is nonzero for every specialization
\beg{equation}\label{rhozerospecialzn}
g=a(s+1)=2ad, \text{ and }m=(a+1)s=(2d-1)(a+1)
\end{equation}
where $a$ and $d$ are positive integers, $a \geq 2$. (Note that the equations \eqref{rhozerospecialzn} encode the fact that $\rho=0$. The possibility that $a=1$ is precluded because in that case the corresponding series $g^{2d-1}_m$ are canonical, and do not determine a divisor in $\mc{G}^{2d-1}_m$, essentially because every canonical curve that admits a $(d-2)$-secant plane admits a one-parameter family of such planes.) Moreover, by Theorem~\ref{hypergeom}, we have
\beg{equation}\label{Pcad}
P_c=P_c(a,d)= -\fr{(2ad)!(2ad-2d+a-1)!}{(2ad-2d)!d!(2ad-d+a-1)!}{ }_3F_2 \biggl[\beg{array}{ccc} -\fr{a}{2}+\fr{1}{2}, & -\fr{a}{2}+1, &-d \\ ad+\fr{1}{2}-d, & ad+1-d \end{array}\biggl| 1\biggr].
\end{equation}
Using \eqref{Pcad}, it is not hard to check that 
\beg{equation}\label{PasQ}
-P_c(a,d)= \fr{(2ad)!}{(2ad-d+a-1)! d!} Q(a,d)
\end{equation}
where 
\beg{equation}\label{Qsum}
Q(a,d)= \sum_{i=0}^{\lfl \fr{a-1}{2} \rfl} (-1)^i \fr{((2a-2)d+a-1)!}{((2a-2)d+2i)!} \cdot \fr{d!}{(d-i)!} \cdot \fr{(a-1)!}{(a-1-2i)!} \cdot \fr{1}{i!}.
\end{equation}
Whenever $a \geq 2$ and $d \geq 1$, the $i$th summand in the sum \eqref{Qsum} has larger absolute value than the $(i+1)$th summand; consequently, $P_c(a,d)$ is negative for all $a \geq 2$ and $d \geq 1$. Nonemptiness follows immediately.
\end{proof}

\subsection{Slopes of secant plane divisors with $r=1$}\label{sloper=1}
Note that for any particular choice of $a \geq 2$, the $i$th summand in \eqref{Qsum} is a polynomial of degree $(a-1-i)$. It follows that
{\small
\[
\beg{split}
P_c(a,d)= -\fr{(2ad)!}{(2ad-d+a-1)! d!} \biggl[\fr{((2a-2)d+a-1)!}{((2a-2)d)!} - \fr{((2a-2)d+a-1)!d}{((2a-2)d+2)!} \cdot (a-1)(a-2)+ O(d^{a-3}) \biggr].
\end{split}
\]
}
Similarly, we have
{\small
\[
\beg{split}
P_{\al}(a,d)&=\fr{(2ad)!}{2(2ad-d+a-1)! d!} \sum_{i=0}^{\lfl \fr{a-1}{2} \rfl} (-1)^i \fr{((2a-2)d+a-1)!}{((2a-2)d+2i)!} \cdot \fr{d!}{(d-i)!} \cdot \fr{(a-1)!}{(a-1-2i)!} \cdot \fr{1}{i!} \\
&- \fr{(2ad-1)!}{2(2ad-d+a-1)! d!} \sum_{i=0}^{\lfl \fr{a}{2} \rfl} (-1)^i \fr{((2a-2)d+a-1)!}{((2a-2)d+2i-1)!} \cdot \fr{d!}{(d-i)!} \cdot \fr{a!}{(a-2i)!} \cdot \fr{1}{i!} \\
&=\fr{(2ad)!}{2(2ad-d+a-1)! d!} \biggl[\fr{((2a-2)d+a-1)!}{((2a-2)d)!} - \fr{((2a-2)d+a-1)!d}{((2a-2)d+2)!} \cdot (a-1)(a-2)+ O(d^{a-3}) \biggr] \\
&-\fr{(2ad-1)!}{2(2ad-d+a-1)! d!} \biggl[\fr{((2a-2)d+a-1)!}{((2a-2)d-1)!} - \fr{((2a-2)d+a-1)!d}{((2a-2)d+1)!} \cdot a(a-1)+ O(d^{a-2}) \biggr],
\end{split}
\]
\[
\beg{split}
P_{\be}(a,d)&= \fr{2(2ad-2)!}{(2ad-d+a-2)!(d-1)!} \sum_{i=0}^{\lfl \fr{a-1}{2} \rfl} (-1)^i \fr{((2a-2)d+a-1)!}{((2a-2)d+2i)!} \cdot \fr{(d-1)!}{(d-1-i)!} \cdot \fr{(a-1)!}{(a-1-2i)!} \cdot \fr{1}{i!} \\
&- \fr{2(2ad-1)!}{(2ad-d+a-1)!(d-1)!} \sum_{i=0}^{\lfl \fr{a-1}{2} \rfl} (-1)^i \fr{((2a-2)d+a)!}{((2a-2)d+2i+1)!} \cdot \fr{(d-1)!}{(d-1-i)!} \cdot \fr{(a-1)!}{(a-1-2i)!} \cdot \fr{1}{i!} \\
&=\fr{2(2ad-2)!}{(2ad-d+a-2)!(d-1)!} \biggl[\fr{((2a-2)d+a-1)!}{((2a-2)d)!} - \fr{((2a-2)d+a-1)!(d-1)}{((2a-2)d+2)!} \cdot (a-1)(a-2)+ O(d^{a-3}) \biggr] \\
&-\fr{2(2ad-1)!}{(2ad-d+a-1)!(d-1)!} \biggl[\fr{((2a-2)d+a)!}{((2a-2)d+1)!}- \fr{((2a-2)d+a)!(d-1)}{((2a-2)d+3)!} \cdot (a-1)(a-2)+ O(d^{a-3}) \biggr],
\end{split}
\]
\[
\beg{split}
P_{\ga}(a,d)&= \fr{8(2ad-5)!}{3(2ad-d+a-3)!(d-3)!} \sum_{i=0}^{\lfl \fr{a-1}{2} \rfl} (-1)^i \fr{((2a-2)d+a)!}{((2a-2)d+2i+1)!} \cdot \fr{(d-3)!}{(d-3-i)!} \cdot \fr{(a-1)!}{(a-1-2i)!} \cdot \fr{1}{i!}\\
&-\fr{7(2ad-2)!}{12(2ad-d+a-1)!(d-1)!} \sum_{i=0}^{\lfl \fr{a}{2} \rfl} (-1)^i \fr{((2a-2)d+a)!}{((2a-2)d+2i)!} \cdot \fr{(d-1)!}{(d-1-i)!} \cdot \fr{a!}{(a-2i)!} \cdot \fr{1}{i!}\\
&+\fr{3(2ad-5)!}{(2ad-d+a-2)!(d-2)!} \sum_{i=0}^{\lfl \fr{a+1}{2} \rfl} (-1)^i \fr{((2a-2)d+a)!}{((2a-2)d+2i-1)!} \cdot \fr{(d-2)!}{(d-2-i)!} \cdot \fr{(a+1)!}{(a+1-2i)!} \cdot \fr{1}{i!}\\
&+\fr{7(2ad-5)!}{12(2ad-d+a-1)!(d-1)!} \sum_{i=0}^{\lfl \fr{a+3}{2} \rfl} (-1)^i \fr{((2a-2)d+a)!}{((2a-2)d+2i-3)!} \cdot \fr{(d-1)!}{(d-1-i)!} \cdot \fr{(a+3)!}{(a+3-2i)!} \cdot \fr{1}{i!}\\
&=\fr{8(2ad-5)!}{3(2ad-d+a-3)!(d-3)!} \biggl[ \fr{((2a-2)d+a)!}{((2a-2)d+1)!}- \fr{((2a-2)d+a)!(d-3)}{((2a-2)d+3)!} \cdot (a-1)(a-2)+ O(d^{a-3}) \biggr] \\
&-\fr{7(2ad-2)!}{12(2ad-d+a-1)!(d-1)!} \biggl[\fr{((2a-2)d+a)!}{((2a-2)d)!}- \fr{((2a-2)d+a)!(d-1)}{((2a-2)d+2)!} \cdot a(a-1)+ O(d^{a-2}) \biggr] \\
&+\fr{3(2ad-5)!}{(2ad-d+a-2)!(d-2)!} \biggl[ \fr{((2a-2)d+a)!}{((2a-2)d-1)!}- \fr{((2a-2)d+a)!(d-2)}{((2a-2)d+1)!} \cdot (a+1)a+ O(d^{a-1}) \biggr] \\
&+\fr{7(2ad-5)!}{12(2ad-d+a-1)!(d-1)!} \biggl[ \fr{((2a-2)d+a)!}{((2a-2)d-3)!}- \fr{((2a-2)d+a)!(d-1)}{((2a-2)d-1)!} \cdot (a+3)(a+2)+ O(d^{a+1}) \biggr],
\end{split}
\]
and
\[
\beg{split}
P_{\de_0}(a,d)&= \fr{8(2ad-5)!}{3(2ad-d+a-3)!(d-3)!} \sum_{i=0}^{\lfl \fr{a-1}{2} \rfl} (-1)^i \fr{((2a-2)d+a)!}{((2a-2)d+2i+1)!} \cdot \fr{(d-3)!}{(d-3-i)!} \cdot \fr{(a-1)!}{(a-1-2i)!} \cdot \fr{1}{i!}\\
&-\fr{(2ad-2)!}{12(2ad-d+a-1)!(d-1)!} \sum_{i=0}^{\lfl \fr{a}{2} \rfl} (-1)^i \fr{((2a-2)d+a)!}{((2a-2)d+2i)!} \cdot \fr{(d-1)!}{(d-1-i)!} \cdot \fr{a!}{(a-2i)!} \cdot \fr{1}{i!}\\
&+\fr{(2ad-5)!}{(2ad-d+a-2)!(d-2)!} \sum_{i=0}^{\lfl \fr{a+1}{2} \rfl} (-1)^i \fr{((2a-2)d+a)!}{((2a-2)d+2i-1)!} \cdot \fr{(d-2)!}{(d-2-i)!} \cdot \fr{(a+1)!}{(a+1-2i)!} \cdot \fr{1}{i!}\\
&+\fr{(2ad-5)!}{12(2ad-d+a-1)!(d-1)!} \sum_{i=0}^{\lfl \fr{a+3}{2} \rfl} (-1)^i \fr{((2a-2)d+a)!}{((2a-2)d+2i-3)!} \cdot \fr{(d-1)!}{(d-1-i)!} \cdot \fr{(a+3)!}{(a+3-2i)!} \cdot \fr{1}{i!}\\
&=\fr{8(2ad-5)!}{3(2ad-d+a-3)!(d-3)!} \biggl[ \fr{((2a-2)d+a)!}{((2a-2)d+1)!}- \fr{((2a-2)d+a)!(d-3)}{((2a-2)d+3)!} \cdot (a-1)(a-2)+ O(d^{a-3}) \biggr] \\
&-\fr{(2ad-2)!}{12(2ad-d+a-1)!(d-1)!} \biggl[\fr{((2a-2)d+a)!}{((2a-2)d)!}- \fr{((2a-2)d+a)!(d-1)}{((2a-2)d+2)!} \cdot a(a-1)+ O(d^{a-2}) \biggr] \\
&+\fr{(2ad-5)!}{(2ad-d+a-2)!(d-2)!} \biggl[ \fr{((2a-2)d+a)!}{((2a-2)d-1)!}- \fr{((2a-2)d+a)!(d-2)}{((2a-2)d+1)!} \cdot (a+1)a+ O(d^{a-1}) \biggr] \\
&+\fr{(2ad-5)!}{12(2ad-d+a-1)!(d-1)!} \biggl[ \fr{((2a-2)d+a)!}{((2a-2)d-3)!}- \fr{((2a-2)d+a)!(d-1)}{((2a-2)d-1)!} \cdot (a+3)(a+2)+ O(d^{a+1}) \biggr].
\end{split}
\]
}

On the other hand, when $r=1$, Khosla's formulas \eqref{Gysin} imply that
\beg{equation}\label{GysinR2}
\beg{split}
\eta_* \al&= -\fr{N(2d-1)(a+1)[(2a^2-2a)d^2 +(a^2+a-8)d+ (4a+2)]}{(2ad-1)(ad-1)} \la\\
&+ \fr{N(2d-1)(a+1)[(2a^2-2a)d^2+ (a^2-4)d+ (2a+1)]}{6(2ad-1)(ad-1)} \de_0, \\
\eta_* \be&= \fr{6N(2d-1)(a+1)}{2ad-1} \la- \fr{N(2d-1)(a+1)}{2(2ad-1)} \de_0, \text{ and } \\
\eta_* c &= -\fr{N(2d-1)[(2a^3-2a)d^3+ (a^3+6a^2-a-8)d^2+ (3a^2+2a-4)d+ a]}{(2d+a)(2ad-1)(ad-1)} \la \\ &+\fr{N(2d-1)d[(2a^3-2a)d^2+ (a^3+4a^2-a-4)d+ (2a^2-2)]}{6(2d+a)(2ad-1)(ad-1)} \de_0.
\end{split}
\end{equation}

Using our hypergeometric formulas for tautological coefficients in tandem with the pushforward formulas \eqref{GysinR2} and \eqref{SecMg}, we may write down the ``virtual slopes" $\fr{b_{\la}}{b_0}$ of secant plane divisors with $r=1$ {\it for any particular value of $a$.} In the following table, we record the virtual slopes corresponding to $2 \leq a \leq 5$.

\medskip
\beg{tabular}{| r | r | r | r | r | r | r | r |} \hline
$a$ & $\fr{b_{\la}}{b_0}$ \\
\hline \hline
2 & $\fr{2(96d^4+80d^3-110d^2-62d+5)}{32d^4+8d^3-30d^2-8d+1}$ \\
\hline
3 & $\fr{3(9216 d^6  + 15552 d^5  + 5240 d^4  - 6372 d^3  - 5218 d^2  - 1067 d + 69)}{4608 d^6  + 6048 d^5  + 772 d^4  - 2780 d^3  - 1609 d^2  - 205 d + 21}$ \\
\hline
4 & $\fr{2 (25920 d^7  + 45360 d^6  + 24387 d^5  - 6006 d^4  - 12143 d^3  - 5213 d^2  - 790 d + 38)}{8640 d^7  + 12744 d^6  + 4853 d^5  - 2585 d^4  - 3032 d^3  - 1041 d^2  - 105 d + 8}$ \\
\hline
5 & $\fr{2 (9830400 d^8  + 18595840 d^7  + 12571776 d^6  + 958200 d^5- 3620196 d^4  - 2433066 d^3  - 734307 d^2  - 89401 d + 3285)}{3276800 d^8+ 5488640 d^7  + 3012992 d^6  - 174328 d^5  - 1038520 d^4  - 575170 d^3  - 145032 d^2- 12207 d + 720}$ \\
\hline
\end{tabular}

\medskip
\medskip
Likewise, our formulas readily yield asymptotics in $d$ for the virtual slopes of secant plane divisors with $r=1$. Namely, we find that when $r=1$, any secant plane divisor on $\ov{\mc{M}}_g= \ov{\mc{M}}_{2ad}$ has virtual slope equal to
\[
\fr{b_{\la}}{b_0}= \fr{6S_1d+ S_2+ O(d^{-1})}{S_1d+S_3+O(d^{-1})}
\]
where
\[
\beg{split}
S_1&=256a^{10}-1024a^9+1280a^8-1280a^6+1024a^5-256a^4, \\
S_2&=384a^{10}+384a^9-13824a^7+768a^8+26496a^6-18048a^5+3072a^4+768a^3, \text{ and }\\
S_3&=64a^{10}-192a^9-2944a^7+1024a^8+3136a^6-448a^5-1152a^4+512a^3.
\end{split}
\]
In particular, we obtain the following result, which describes the $d$-asymptotics of secant plane divisors when $r=1$.
\beg{thm}
For every choice of positive integers $(a,d)$, the difference between the $(a,d)$th secant plane divisor's (virtual) slope and that of the Brill-Noether divisor on $\ov{\mc{M}}_{2ad}$ is equal to
\[
\fr{b_{\la}}{b_0}- 6- \fr{12}{2ad+1}= \fr{3}{ad(a+1)}+ O(d^{-2})= \fr{6}{(a+1)g}+ O(g^{-2}).
\]
\end{thm}

\section{Boundary coefficients of secant plane divisors on $\ov{\mc{M}}_g$}\label{bdrycoeffs}
Much as in the preceding section, write the class of $\mbox{Sec}$ as an expansion in terms of standard divisor classes on $\ov{\mc{M}}_g$:
\[
\mbox{Sec}= b_{\la} \la- b_0 \de_0- \sum_{i=1}^{\lfl \fr{g}{2} \rfl} b_i \de_i.
\]
In this section, we will determine $b_1$ and $b_2$.

\subsection{Determination of $b_1$}\label{b1}
Consider the curve $Y \approx \mb{P}^1 \hra \ov{\mc{M}}_g$ given by attaching a general pencil of plane cubics to a general genus-$g$ flag curve $Y$ at a general point of $Y$. By the same argument used to prove \cite[Thm 1]{C}, we see that $Y$ avoids every secant plane divisor. On the other hand, it is well-known (see, e.g., \cite[p. 174]{HM}) that
\[
Y \cdot \la=1, Y \cdot \de_0=12, Y \cdot \de_1=-1, \text{ and } Y \cdot \de_i=0 \text{ for all }i \geq 2.
\]
It follows that 
\[
b_{\la} - 12 b_0+b_1= 0.
\]

\subsection{Determination of $b_2$}\label{b2}
Given any integer $\al \geq 2$, let
\[
j_{\al}: \ov{\mc{M}}_{\al,1} \ra \ov{\mc{M}}_g
\]
denote the map defined by attaching a fixed flag curve $C$ of genus $(g-\al)$ at a fixed general point  of $C$ to any genus-$\al$ curve $Y$ with a marked point. Much as in \cite[proof of Thm 1.1]{FP}, whose argument we follow, we have the following result.
\beg{thm}\label{delta2coeff}
If $\al=2$, then $j_{\al}^* \mbox{Sec}$ is supported on the Weierstrass locus.
\end{thm}
Recall that the Weierstrass locus comprises curves marked along Weierstrass points, and has class
\[
\ov{\mc{W}}= -\la+ \fr{g(g+1)}{2} \de_0- \sum_{i=1}^{g-1} \binom{g-i+1}{2} \de_i
\]
according to \cite{Cu}. 

\beg{proof}
Assume, for the sake of argument, that $j_{\al}^* \mbox{Sec}$ is {\it not} supported on the Weierstrass locus; this means exactly that some curve $C \cup_p Y$, where $p$ is not a Weierstrass point of $Y$, carries a pair of limit linear series $(g^s_m,g^{s-d+r}_m)$ satisfying \eqref{seriesinseries}. Moreover, by additivity of the generalized Brill--Noether number, we have
\beg{equation}\label{add1}
\rho(2,s,m; r(Y,p))+ \rho(g-2,s,m; r(C,p))= \rho(g,s,m)=0
\end{equation}
where $r(Y,p)$ and $r(C,p)$ denote the total ramification of the $g^{s-d+r}_m$ along $Y$ and $C$, respectively. Since $(C,p)$ is Brill--Noether general by assumption, and $(Y,p)$ is Brill--Noether general whenever $p$ is not a Weierstrass point of $Y$, \eqref{add1} forces
\beg{equation}\label{add2}
\rho(2,s,m; r(Y,p))= \rho(g-2,s,m; r(C,p))=0.
\end{equation}
Because $p$ is not a Weierstrass point of $Y$, we now deduce that the vanishing sequence at $p$ of the aspect of the $g^s_m$ along $Y$ is either
\[
\beg{split}
&a(V_Y,p)= (m-s-2, m-s-1, \dots, m-3,m), \text{ or } \\
&a(V_Y,p)= (m-s-2, m-s-1, \dots, m-4,m-2,m-1).
\end{split}
\]
Now assume that $a$ base points of the included series $g^{s-d+r}_m$ lie along $Y$. Thus, $(d-a)$ base points lie along $C$, which in turn forces $(d-a)(s-d+r)$ shifts of vanishing order indices of the $g^{s-d+r}_m$ along $C$, as shown in \cite[Lem. 1]{C} (from which the terminology of ``shifting" is borrowed as well). On the other hand, from standard theory of limit linear series \cite{EH3}, we have
\beg{equation}\label{limitlinearfact}
a_i(V_Y,p)+ a_{s-i}(V_C,p) \geq m
\end{equation}
for all $0 \leq i \leq s$. Via \eqref{limitlinearfact}, the $(d-a)(s-d+r)$ shifts of vanishing order indices of the $g^{s-d+r}_m$ along $C$ impose at least $(s-d+r)(d-a)$ degrees of ramification along $Y$ away from $p$. It follows that the total ramification of the $g^{s-d+r}_{m-a}$ along $Y$ obtained by removing the $a$ base points from our $g^{s-d+r}_m$ is at least 
\[
r=(s-d+r+1)(m-s-2)+ (s-d+r)(d-a).
\]
But an easy calculation yields
\[
\beg{split}
\rho(2,s-d+r,m-a)- r &= 2+ (r-d-1)- (a-d+r), \text{ since } \mu(d,s,r)=-1 \\
&= 1-a.
\end{split}
\]
Because $(Y,p)$ is Brill--Noether general, it follows that $a \leq 1$. Clearly, the case $a=0$ is impossible, since this would imply that every base point of the $g^{s-d+r}_m$ lies along $C$, and, therefore, that $C$ admits a $d$-secant $(d-r-1)$-plane. On the other hand, the case $a=1$ is also precluded, because in that situation \eqref{limitlinearfact} forces the top two vanishing orders at $p$ of the $g^{s-d+r}_m$ along $Y$ to be maximal. This is clearly impossible when 
\[
a(V_Y,p)= (m-s-2, m-s-1, \dots, m-3,m);
\]
for, if there is a base point along $Y$, then the order to which the $g^{s-d+r}_m$ along $Y$ vanishes at $p$ must be less than $m$. Similarly, if 
\[
a(V_Y,p)= (m-s-2, m-s-1, \dots, m-4,m-2,m-1),
\]
then subtracting the base point from the $g^{s-d+r}_m$ along $Y$ yields a $g^{s-d+r}_{m-1}$ containing a subpencil $\Ga$ of sections vanishing to orders $(m-2)$ and $(m-1)$, respectively. Subtracting $(m-2)$ base points from $\Ga$ yields a $g^1_1$ along the genus-2 curve $Y$, which is absurd.
\end{proof}

As explained in \cite[Thm. 6.65]{HM}, Theorem~\ref{delta2coeff} implies that
\[
b_2= \fr{5}{2}b_1- \fr{b_{\la}}{2}.
\]

It is natural to ask whether a suitable modification of the argument used to prove Theorem~\ref{delta2coeff} may be used to determine any of the remaining boundary divisor coefficients $b_i, i \geq 3$. The basic issue here is whether there exist reducible curves $C \cup_p Y$ of genus $g$, where $(C,p)$ is a general pointed flag curve and $(Y,p)$ is a Brill--Noether general pointed curve, that admit exceptional secant planes. In other words, one would like to see how the pullbacks $j_i^* \mbox{Sec}$ relate to the Brill--Noether divisors on $\ov{\mc{M}}_{i,1}$ studied in \cite{EH4}.

{\fl \bf Question:} For which values of $i \geq 3$ and $(d,r)$ is it the case that the pullback of $\mbox{Sec}$ along $j_i$ is supported along a union of Brill--Noether divisors on $\ov{\mc{M}}_{i,1}$?

As noted in \cite{FP}, if $j_i^* \mbox{Sec}$ is supported along a union of Brill--Noether divisors on $\ov{\mc{M}}_{i,1}$ for every $i \leq j$, then every boundary divisor coefficient $b_i, i \leq j$, because, as is shown in \cite{EH4}, the class of every Brill--Noether divisor on $\ov{\mc{M}}_{i,1}$ is a linear combination of $\ov{\mc{W}}$ and the pullback of the Brill--Noether divisor on $\ov{\mc{M}}_i$.

{\fl \bf \small Laboratoire de Math\'ematiques Jean Leray \\
Unit\'e mixte de recherche 6629 du CNRS\\
Universit\'e de Nantes\\
2 rue de la Houssini\`ere, BP 92208 \\
44322 Nantes Cedex 3, FRANCE \\

{\it Email address}: \url{cotteril@math.harvard.edu} \\
{\it Web page}: \url{www.mast.queensu.ca/~cotteril}}

\end{document}